\documentclass[12pt,reqno]{amsart}

\usepackage{amsmath,amsfonts,amsthm,amssymb,color,hyperref}
\usepackage{lmodern}

\usepackage{amsthm}
\usepackage{amssymb} 
\usepackage{cmmib57} 
\usepackage{exscale} 
\usepackage{amscd} 
\usepackage{latexsym} 

\usepackage{graphicx}
\usepackage{pdfsync}
\usepackage{wasysym}
\usepackage[catalan,spanish,english]{babel} 
\usepackage{textcomp} 
\usepackage{bm}  
\usepackage{colortbl}
\usepackage{fancyhdr}
\usepackage{cancel}
\usepackage{scalerel,stackengine}
\usepackage{cite}

\usepackage{float}
\usepackage[utf8]{inputenc}
\usepackage{mathrsfs}

\usepackage[T1]{fontenc}

\usepackage{subfigure} 

\makeatletter
     \def\section{\@startsection{section}{1}%
     \z@{.7\linespacing\@plus\linespacing}{.5\linespacing}%
     {\bfseries
     \centering
     }}
     \def\@secnumfont{\bfseries}
     \makeatother

\newcommand{\IN}{{\mathbb N}}

\newcommand{\beq}{\begin{equation}}
\newcommand{\eeq}{\end{equation}}
\newcommand{\bal}{\begin{align}}
\newcommand{\eal}{\end{align}}
\newcommand{\beqn}{\begin{equation*}}
\newcommand{\eeqn}{\end{equation*}}
\newcommand{\baln}{\begin{align*}}
\newcommand{\ealn}{\end{align*}}

\newcommand{\E}{\mathbb{E}}

 \allowdisplaybreaks

 \DeclareGraphicsExtensions{.pdf,.png,.jpg,.tiff}

\restylefloat{figure}

\title[Modelling mortality rates with a geometric-type fOU]{Modelling Italian mortality rates with a geometric-type fractional Ornstein-Uhlenbeck process.}

\author[{\small F. Delgado-Vences}]{Francisco Delgado-Vences}
\address{Conacyt Research Fellow - Universidad Nacional Aut\'onoma de M\'exico. Instituto de Matem\'aticas, Oaxaca, M\'exico}
\email{delgado@im.unam.mx}

\author[{\small A. Ornelas}]{Arelly Ornelas} 
\address{Conacyt Research Fellow - Instituto Politecnico Nacional- CICIMAR, La Paz, M\'exico}
\email{arelly.ornelas@conacyt.mx}

\newtheorem{theorem}{Theorem}[section]

\theoremstyle{definition}

\theoremstyle{remark}

\numberwithin{equation}{section}
\begin{document}

\maketitle


\allowdisplaybreaks


\begin{abstract}
We propose to model mortality hazard rates for human population using the exponential of the solution of an stochastic 
differential equation (SDE). The noise in the SDE is a fractional Brownian motion. We will use the 
 well-known fractional Ornstein-Uhlenbeck process. Using the Hurst parameter we showed that
  mortality rates exhibit long-term memory. 

The proposed model is a generalization of the model introduced by 
\cite{gi-or-be}, where they used an SDE driven with a Brownian motion. 
We tested our model with the Italian population between 
the years 1950 to 2004.

\end{abstract}

\section{Introduction}
\label{intro}

Future planning in the demographic, economic and actuarial areas is crucial, for instance: good planning in 
social programs, government budgets, actuarial reserves, cost of insurance and pensions, etc., 
depends on the use of a good method to forecast. However,  constant changes in technology, 
lifestyle, climate change, migration, to name a few, make  predicting  a non easy task. These changes
have given rise to the need for further research in the  field of risk management, insofar as life expectancy affects mortality
forecasting, longevity risk, reserve calculations, annuities, pension plan design and premiums for life products.
Mortality  impacts directly in  cash costs and therefore need a good future projection. \\

Several models have been proposed to describe mortality. Pitacco et al. (2009) provides an interesting review of the 
historical background to early mortality tables. Some authors in their attempts to calibrate models to mortality rates
make the assumption of perfect correlation across generations. However, one could think that correlation among close generations is
high but not perfect. This hypothesis is 
our start point for this paper, we pretend to model the mortality hazard rates for the human population considering a 
model that include the high correlation among generations.
In \cite{je-lu-vi}, they study a cohort-based model which uses the imperfect correlation of mortality intensity
across generations. They implemented  it on UK data for the period 1900-2008.\\

 In this section we present a brief review of the Milevesky-Promislow model and  discuss about long-range dependence. 
 Milevesky-Promislow model the mortality hazard rate using the 
 exponential of an Stochastic Differential Equation (see \cite{mi-pr}).  \\

Stochastic Differential Equations (SDEs) provide a powerful and flexible probabilistic structure for modeling phenomena in a 
multitude of disciplines: finance, biology, molecular dynamics, chemistry, survival analysis, epidemiology, just to name a few
examples, the reader could see \cite{ok} or \cite{kl-pl} for further information. The fractional Ornstein-Uhlenbeck process 
(FOU) process we are interested is an example of an SDE.\\

In \cite{gi-or-be}, the authors present a generalization of the Milevesky-Promislow model that covers the case with diffusion 
coefficient as a function of $t$ instead of  a constant like in \cite{mi-pr}. Also, they consider a type of Autoregressive model 
for the logarithm of the hazard rate\footnote{see equations (13) or (16) in \cite{gi-or-be}}. It is known that all stationary invertible 
ARMA(p,q) processes are short memory processes\footnote{see for instance page 739 in \cite{be-ye-fe-gh} or section 5.2 
in \cite{sh-st}}.
 In particular the authors in  \cite{gi-or-be} used a modified AR(1) process, that continuous being a 
short memory process. In our model we consider the case  with a long memory process. For another generalization of the 
Milevesky-Promislow model see also \cite{ro-le}. \\

In the other hand we are interested in the long range dependence. 
Hurst \cite{hu} observed a phenomenon which is invariant to changes in scale, when
he was studying the river water levels along the Nile. Another scale-invariant phenomenon was
observed in studies of problems connected with traffic patterns of packet
flows in high-speed data networks such as the Internet (see preface of \cite{ra}).\\

Long Range dependence (LRD) ( also known as long memory, strong dependence or persistence)
denotes the property of a stochastic process or a time series to exhibit persistent behavior across long time periods 
(see \cite{ra}). LRD is used in areas such as finance, econometrics, 
Internet, hydrology, climate studies, linguistics, geophysics or DNA sequencing among others. \\

One example of time-continuous long-range dependence model was introduced by
Mandelbrot and van Ness (see \cite{ma-va}), they define the term
fractional Brownian motion (fBm) for a Gaussian process with a specific covariance
structure and studied its properties. This process is a generalization of classical
Brownian motion also known as the Wiener process (\cite{ra}). The self-similarity and 
long-range dependence properties make the fractional Brownian motion 
suitable to model driving noises in different applications such as hydrology, Short-term 
Rainfall Prediction, finance, etc.. \\

The fBm is a particular case of the Fractional Gaussian noises (FGN), that are a type of stochastic 
processes. The FGN have the property of self-similarity and is
used to model persistent dependency in time series. 
The autocovariance function of FGN is characterized by the Hurst exponent (H parameter). 
fBM and fOU are examples of fractional Gaussian noise. \\

We now discuss the model of Milevsky-Promislow Model (see \cite{mi-pr} for the original paper 
or \cite{gi-or-be} 
for a recent generalization). \\

Set the survival probability of an individual aged $x$ in the period $[t,T]$ as 
\begin{align}\label{E-surv}
 S(t,T):=\E\Big[\exp\Big(-\int_t^T h_x(u) du\Big)\Big|\mathcal{F}_t \Big],
\end{align}
where $\{\mathcal{F}_t\}_{t\ge 0} $ is a filtration which represent the information until time $t$ and $h_x(t)$ is the 
stochastic force of mortality or hazard rate. According with the Milevsky-Promislow Model $h_x(t)$ is given by
\begin{align}
 h(t)&=h_0\exp(\alpha_0t+\alpha_1Y_t),\label{mod-1}
\end{align}
where  $h_0,\alpha_1,\alpha_2> 0$. The process $Y_t$ satisfies the SDE:
\begin{align}
 dY_t&=-\lambda Y_tdt+\sigma dB_t, \label{mod-3}
\end{align}
where $B_t$ is a Brownian motion, $Y_0=0$ and $\sigma,\lambda> 0$.\\

In this paper we assume, according with the Milevsky-Promislow Model, that $h_t$ is given by 
\eqref{mod-1}. However $Y_t$ satisfies an SDE whose solution is the fractional Ornstein-Uhlenbeck process. Indeed,
we will assume that the stochastic process $Y_t$ satisfies the following stochastic differential
equation (SDE):
\begin{align}
 dY_t^H&=-\lambda Y_t^Hdt+\sigma dB_t^H,  \label{mod-2}
\end{align}
Where $B_t^H $ is a fBM with Hurst parameter $1/2 \le H< 1$, $Y_0^H=0$, and $\sigma,\lambda> 0$.  
$Y_t^H$ is called a fractional Ornstein-Uhlenbeck process (fOU) and has been studied exhaustively 
in the last decades. \\

A priori, we will assume that $\alpha_1=1$. Once we estimated the parameter $H$, we will 
adjust $\alpha_1=T^{-H}$, which is devoted to control the variance of the process $Y_t$. 
We will use $Y_t$ to denote $Y_t^H$.\\

Since stochastic mortality rate models take into account long time phenomena, we  
suggest a generalization of Milevesky-Promislow model, given by the equations 
\eqref{mod-1} and \eqref{mod-2}, which cover the case of long range dependence of
the data.\\

The main difference between this model and the original presented 
in \cite{mi-pr} is that they consider the driving noise in the SDE as a standard
Brownian motion instead of the fractional Brownian motion as in our model.\\

The paper is organized as follows. In section \ref{fgn} we present a brief review of 
the fBM and fOU and we discuss some important properties of the fOU. Section \ref{esti}
is dedicated to describe parameters estimation. 
In first place we discuss the estimation of parameter $\alpha_0$. In second place we review several
methods for estimating  Hurst parameter 
$H$. Finally we present some results devoted to estimate the parameters 
$\sigma,\lambda$ in the SDE. The results of the 
model are discussed in section \ref{re-fou}.\\

\section{On the fractional Gaussian noise}
\label{fgn}

\subsection{Fractional Brownian motion}
A stochastic process is a collection of random variables indexed by time. It can be at discrete  or continuous time. 
A discrete time stochastic process $X = \{X_n, n = 0, 1, 2,\ldots\}$ is a countable collection of random variables 
indexed by the non-negative integers. A continuous time stochastic process $X = \{X_t , 0\le t <\infty\}$ is an 
uncountable collection of random variables indexed by the non-negative real numbers. The most notable example of 
stochastic process at continuous time is the so-called the Brownian motion (BM), which is a collection of random 
variables $\{B_t\}$ such that
\begin{itemize}
 \item $B_0=0$ almost surely.
 \item For $0\le s < t < \infty$,\quad $B_t-B_s \sim N(0, t-s)$,
 \item For $0\le s < t < \infty$,\quad $B_t-B_s$ is independent of $B_s$,
\item The trajectories $t\mapsto B_t$ are continuous.
 \end{itemize}
The Brownian motion is one of the most important stochastic process and is an example of a Gaussian process.  \\

We consider a generalization of the Brownian motion. Let $\{B_t^H,t\ge 0\}$, with $H\in (0,1)$, be a Gaussian process with zero-mean 
and covariance function given by
\begin{align}
 R_H(t,s):=\E(B_s^HB_t^H)= \tfrac{1}{2}\big(t^{2H}+s^{2H}-|t-s|^{2H} \big).\label{s1.1}
\end{align}
This stochastic process  is called a {\it fractional Brownian motion} (fBm) and was introduced by Kolmogorov \cite{ko} 
and studied by Mandelbrot and Van Ness in \cite{ma-va}. The parameter $H$ is called Hurst index because of  
the statistical analysis developed by the climatologist Hurst \cite{hu}.
The fBm is a generalization of Brownian motion without independent increments, also it 
is a continuous-time Gaussian process. \\

The fBm has the following properties:

\begin{enumerate}
 \item {\bf Self-similarity}: The processes $\{a^{-H} B_{at}^H, t \ge 0\}$ and $\{B_t^H,t\ge 0\}$
 have the same probability distribution, for any constant $a > 0$.
\item {\bf Stationary increments}: From \eqref{s1.1} it follows that the increment of the process
in an interval $[s,t]$ has a normal distribution with zero mean and variance equal to:
\begin{align*}
 \E\big((B_t^H-B_s^H)^2 \big)= |t-s|^{2H}.
\end{align*}
\item {\bf Sample-paths are almost nowhere differentiable}. However, almost-all trajectories are
H\"older continuous of any order strictly less than H: for each such trajectory, there exists a
finite constant $C$ such that
\[
 \E\big(|B_t^H-B_s^H|\big) \le C |t-s|^{H-\epsilon}, 
\]
for every $\epsilon > 0$.
\end{enumerate}

For $H =\tfrac{1}{2}$ the covariance can be written as $R_{1/2}(t,s) = \min(s,t)$ and the process $B_t^{1/2}$ is 
the ordinary Brownian motion (Bm). The increments of this process, in disjoint intervals, are independent. However, 
the increments are {\it not} independent for $H \ne\tfrac{1}{2}$.\\

Set  $X_n = B_n^H-B_{n-1}^H$, $n \ge 1$ a  stochastic process.  $\{X_n, n \ge 1\}$ is a Gaussian stationary sequence with unit
variance and covariance function: 
\begin{align*}
 \rho_H (n) &= \frac{1}{2}\Big((n+1)^{2H}+ (n - 1)^{2H}-(2n)^{2H}\Big)\\
 &\approx H (2H - 1)n^{2H-2} \rightarrow 0, \quad \hbox{\rm  when  } n\rightarrow\infty.
\end{align*}
Therefore, 
\begin{itemize}
 \item if $H > \tfrac{1}{2}$ then $\rho_H(n) > 0$ for $n$ large enough and $\sum_{n=1}^\infty \rho_H(n)=\infty$.
This is persistent process with positive correlation. In this case, we say that $X_n$ has 
{\it long-Range dependence} property.
\item If $H < \tfrac{1}{2}$, then $\rho_H(n) < 0$ for $n$ large enough and $\sum_{n=1}^\infty \rho_H(n)<\infty$. This is an
anti-persistent process with  negative correlation.
\end{itemize}
For further information on fBM see  \cite{ra}, \cite{nu} or \cite{mi}.

\subsection{Fractional Ornstein-Uhlenbeck process (fOU)}

The fOU is an SDE driven by a fractional Brownian motion. The same model was used in \cite{mi-pr} 
or \cite{gi-or-be}.\\

As we mentioned before, the survival probability $S(t,T)$ of an individual aged $x$ in the period $[t,T]$,   
is given in the Equation \eqref{E-surv} 
and $h_x(t)$ is the stochastic force of mortality or hazard rate given by the Equation \eqref{mod-1}.
We will assume that $Y_t$ is an stochastic process that satisfies the SDE:
\begin{align}
 dY_t^H&=-\lambda Y_t^Hdt+\sigma dB_t^H, \label{mod2}
\end{align}
where $B_t^H $ is a fBM with Hurst parameter $1/2 \le H< 1$,  $Y_0=0$, and $\sigma,\lambda> 0$.  
This SDE is the fractional Ornstein-Uhlenbeck process.\\

There are substantial differences in trying to solve Equation \eqref{mod2} with respect to the method use
in \cite{mi-pr} or \cite{gi-or-be}. We now discuss some of these differences. First we interpret the SDE \eqref{mod2} as 
\begin{align}
 Y_t^H&=-\lambda\int_0^t Y_s^Hds+\sigma B_t^H.\label{mod3}
\end{align}
Notice that the equation above does not have an stochastic integral because we are considering the 
case with additive noise.
Nevertheless, it is possible to consider the general case with multiplicative noise in which case 
it is necessary to define an stochastic integral with respect to fractional Brownian motion as a pathwise Riemann–Stieltjes 
integral (see, e.g., \cite{yo} for the original definition and \cite{du-no} for advanced results). \\

Coming back to Equation \eqref{mod3}, Cheridito et al in \cite{ch-ka-ma} have introduced the 
fractional Ornstein-Uhlenbeck process (fOU) and they have shown that the process 
\begin{equation}
 Y_t^H=\sigma \int_0^t e^{-\lambda(t-u)} dB_u^H,\label{mod4}
\end{equation}
is the unique a.s. continuous-path process which solves \eqref{mod3} (see also Theorem 1.24 in \cite{ra}).
The integral in Equation \eqref{mod4} is a 
pathwise Riemann–Stieltjes integral.  The fOU process is neither Markovian nor a semimartingale for $H \in(1/2,1)$ but remains
Gaussian and ergodic.\\

Moreover, when $H \in(1/2,1)$, $Y_t$ even presents the long-range dependence property (see Cheridito 
\cite{ch-ka-ma} or  \cite{ra}).\\

The variance of the fOU process $Y_t$ is given by the following expression (see\cite{ze-ch-ya}):
\begin{align}
 Var(Y_t)= \sigma^2 2H e^{-2\lambda t} \int_0^t s^{2H-1} e^{2\lambda s} ds.\label{var-fou}
\end{align}
Notice that when $H=1/2$ we get 
\begin{align}
 Var(Y_t)= \frac{\sigma^2}{2\lambda}  \big(1-e^{-2\lambda t}\big),
\end{align}
which is the variance of the standard Ornstein-Uhlenbeck process (see for instance \cite{mik} page 143). \\

If we consider the constant $\alpha_1 = T^{-H}$ and the Equation \eqref{var-fou},
the expression for the variance of $\alpha_1 Y_t$ is given by
\begin{align}
  Var(\alpha_1 Y_t)&= \alpha_1^2 Var(Y_t)= \alpha_1^2\sigma^2 2H  \int_0^t s^{2H-1} e^{-2\lambda (t-s)} ds\nonumber
 \\
 &\le \alpha_1^2 \sigma^2 2H  \int_0^t s^{2H-1} ds= \alpha_1^2 \sigma^2 2H  \frac{s^{2H}}{2H}\Big|_{s=0}^t \nonumber
 \\
  &= \alpha_1^2 \sigma^2 t^{2H} = \sigma^2 (t/T)^{2H},\label{var-fou1}
\end{align}
then $Var(\alpha_1 Y_t)\le \sigma^2 $ since $0\le t\le T$, this implies that the variance of $\alpha_1 Y_t$ 
is bounded by a constant that does not depend on time. We will use $\alpha_1$  to control the variance of the process $Y_t$.

\section{Estimation of the parameters}
\label{esti}

In this section we will describe a methodology to estimate the parameters. 
We will use some R-libraries in order to estimate the parameters. \\

We need to estimate $\alpha_0, \alpha_1$ as well as $\sigma,\lambda$ for the SDE model that we described in the previous section. Furthermore,
The Hurst parameter ($H$) involved in the driven fractional Brownian motion will be also estimated, however, this estimation is highly complicated.
To solve this problem we will use the empirical evidence that the Hurst value in the equation \eqref{mod2} 
is preserved, this means that the  value of the  Hurst parameters $H$, in the equation \eqref{mod2}, for the fBm $B_t^H$ and the one for the 
fractional Gaussian noise $Y_t^H$ are the same. Observe that $\alpha_1=T^{-H}$ will be calculated using the estimated Hurst parameter.

\subsection{Estimation of the parameter $\alpha_0$.}

The model is given by the equation \eqref{mod-1}. In order to estimate $\alpha_0$ we will assume that $\alpha_1=1$. Taking $\ln$ we obtain

\begin{align}
 \ln h(t)&=\ln h_0+\alpha_0t+Y_t. \label{mod-ln}
\end{align}

One simple method to estimate the parameter $\alpha_0$ is by minimizing the sum of the square errors. 
Let $T$ be given by
\begin{align*}
 T:= \sum_{t_{initial}}^{t_{final}} \Big( \ln h(t)-\ln h_0-\alpha_0t \Big)^2. 
\end{align*}
Taking derivative of $T$ with respect 
to $\alpha_0$  we get
\begin{align*}
 \frac{\partial T}{\partial \alpha_0}&= -2\sum_{t_{initial}}^{t_{final}} \Big( \ln h(t)-\ln h_0-\alpha_0t \Big) t =0, 
\end{align*}
and from this equation we obtain $\widehat{\alpha_0}$:
\begin{align}
\widehat{\alpha_0} &= \frac{\sum_{t} t\ln h(t) - \ln h(0)\sum_{t} t}{\sum_{t} t^2}. \label{alpha0}
\end{align}

Once we have estimate $\alpha_0$ we proceed to estimate the Hurst parameter, $\sigma$ and $\lambda$.

\subsection{Relation between the Hurst parameter and the H-index in the FOU}

In this subsection we will discuss the procedure we have used to estimate the parameter $H$.\\

We will use the following empirical fact. Suppose that a fOU process is driving with a fBM with a given 
Hurst parameter $H_0$.  Yerlikaya-Okzurt et al (\cite{ye-etal})  have show a relationship between the Hurst parameter $H$ of the fractional 
Brownian motion and the Hurst parameter of the fractional Gaussian noise given by an SDE. In fact, they have found statistical evidence 
that the fOU should have the same value $H_0$ that the fBM (see table 1 Yerlikaya-Okzurt et al). Therefore, at least empirically, 
the value of $H$ is the same. Then, it is possible to choose the same value of the 
parameter $H$ for both processes.  A formal proof of this fact, up to our knowledge, is missed. \\

The subsequent sections are devoted to present several methods to estimate the Hurst parameter for the fBM.

\subsection{Estimation of the self-similarity index $H$ for the fBM}
\label{Desc-Est-H}

The last subsection allows us to estimate the parameter $H$ in one simple way. According to equation \eqref{mod-ln}, the residuals are given
by the expression
\begin{align*}
\hat Y_t= \ln h(t)&-\ln h_0 -\widehat \alpha_0t.
\end{align*}
$\hat Y_t$ is a fractional Gaussian noise, so that we can use it to estimate $ H$. Afterwards we will use  $\hat H$ to approximate 
the Hurst parameters of the fractional Brownian motion $B_t^H$. 
For this purpose we will review some methods to estimate
the parameter $H$.

\subsubsection{R over S Analysis}

Following \cite{we}. The analysis begins  dividing a time series $\{Z_i\}$ of length $L$ into $d$
subseries of length $n$ and denote it by $\{Z_{i,m}\}$, $m = 1,\ldots,d$. Then, for each subseries  $\{Z_{i,m}\}$, $m = 1,\ldots,d$ :
\begin{enumerate}
 \item Find the mean $E_m$ and standard deviation $S_m$.
\item Normalize the data $Z_{i,m}$  by subtracting the sample mean $X_{i,m} =Z_{i,m}-E_m$ for $i=1,\ldots,n$. 
\item Create a cumulative time series $Y_{i,m} = \sum_{j=1}^i X_{j,m}$.
for $i = 1,\ldots,n$ 
\item Find the range $R_m = max\{Y_{1,m},\ldots, Y_{n,m}\} - min\{Y_{1,m},\ldots, Y_{n,m} \}$; 
\item Rescale the range $R_m/S_m$ . 

\item Calculate the mean value of the rescaled range
for all subseries of length n
\[
 (R/S)_n =\frac{1}{d}\sum_{m=1}^d  R_m /S_m.
\]
\end{enumerate}

It can be shown (see \cite{we}) that the $R/S$ statistic asymptotically follows the relation:
\[
(R/S)_n \sim c n^H,
\]
where $c$ is a constant. Thus, the value of $H$ can be obtained by running a simple linear regression over a
sample of increasing time horizons
\[
\log(R/S)_n = \log c + H \log n.
\]

Equivalently, we can plot the $(R/S)_n$ statistics against $n$ on a double-logarithmic paper.
If the returns process is white noise then the plot is roughly a straight line
with slope $0.5$. If the process is persistent then the slope $H$ is greater than $ 0.5$; if it is anti-persistent
then the slope $H$ is less than  $0.5$. The “significance” level of the estimated parameter $H$ is usually chosen to be one over the
square root of sample length, i.e. the standard deviation of a Gaussian white noise.\\

A major drawback of the $R/S$ analysis is  that no asymptotic distribution
theory has been derived for the Hurst parameter $H$ . The only results known are for
the rescaled (but not by standard deviation) range $R_m$ itself, see \cite{lo}.

\subsubsection{Method of rescaled range analysis R/S}
Following \cite{ra}, chapter 9.
This method was suggested by Hurst (1951). The series $\{X_j , 1 \le j \le N - 2\}$ is
divided into $K$ nonoverlapping blocks such that each block contains $M$ elements
where $M$ is the integer part of $N/K$. Let $t_i = M(i - 1)$ and
$$
R(t_i,r) = \max[W(t_1 , 1), \ldots, W(t_i , r)] - \min[W (t_1 , 1),  \ldots , W (t_i , r)],
$$
where
\[
 W (t_i , k) = \sum_{j=0}^{k-1} X_{t_i +j}  - k \Bigg(\frac{1}{r}\sum_{j=0}^{r-1} X_{t_i +j}  \Bigg),\quad k = 1,\ldots , r.
\]

Note that $R(t_i, r) \ge 0$ since $W(t_i , r) = 0$ and the quantity $R(t_i , r)$ can be computed only when $t_i + r \le N$ . Define
\[
 S^2 (t_i , r) = \frac{1}{r} \sum_{j=0}^{r-1} X_{t_i +j}^2  -  \Bigg(\frac{1}{r}\sum_{j=0}^{r-1} X_{t_i +j}  \Bigg)^2.
\]
The ratio $R(t_i , r)/S(t_i , r)$ is called the rescaled adjusted range. It is computed
for a number of values of $r$ where $t_i = M(i - 1)$ is the starting point of the ith
block for $i = 1,\ldots, K$. Observe that, for each value of $r$, we obtain a number
of $R/S$ samples. The number of samples decrease as $r$ increases. However, the
resulting samples are not independent. It is believed that the $R/S$-statistic is
proportional to $r^H$ as $r \rightarrow \infty$ for the fractional Gaussian noise. Assuming this
property, we regress $log(R/S)$ against $log(r)$ to obtain an estimator for $H$.

\subsubsection{FDWhittle Estimator}

Following\cite{pa-etal}.nThe Local Whittle Estimator (LWE) is a semiparametric Hurst parameter estimator based on the periodogram.
It assumes that the spectral density $f(\lambda)$ of the process can be
approximated by the function
\begin{equation}\label{fdw1}
f_{c,H}(\lambda) = c\lambda^{1-2H},
\end{equation}

for frequencies $\lambda$ in a neighborhood of the origin, $c$ is a constant.  The periodogram of a time series 
$\{X_t , 1 \ge t \ge N \}$ is defined by 
$$
 IN (\lambda)  = \frac{1}{2\pi N}\left| \sum_{t=1}^N X_te^{i\lambda t}\right|^2,
 $$

where $i=\sqrt{-1}$. Usually, it is evaluated at the Fourier Frequencies 
$\lambda_{j,N} = \frac{2\pi j}{N}$, $0 \le j \le [N/2]$.  Note that the periodogram is the norm of the Discrete Fourier transform of the time series
(see Section 6.1.2 in \cite{priestley} for instance).\\

 The LWE of the Hurst parameter, $\hat{H}_{LWE}(m)$  is implicitly defined by minimizing
 $$
 \sum_{j=1}^m  log f_{c,H} (\lambda j,N ) + \frac{I_N(\lambda j, N)}{ f_{c,H} (\lambda j,N )},
 $$
 with respect to $c$ and $H$, with $f_{c,H}$ defined in \eqref{fdw1}.

\subsection{Estimation of $\sigma$ and  $\lambda$}

There are several methods to estimate parameters $\sigma$ and  $\lambda$. For 
instance see \cite{ra} or the references in \cite{ne-ti} or in \cite{ku-mi}. In the following section we will do a brief review of 
some of these methods.\\

\subsubsection{Estimation $\sigma$ with quadratic generalized variations method.}
\label{sect-est}

 Brouste and Iacus \cite{br-ia} proposed a consistent and asymptotically Gaussian estimators for the parameters $\sigma,\lambda$ and $H$ of 
 the discretely observed fractional Ornstein-Uhlenbeck process solution of the stochastic differential equation. 
 There is a restriction on  the estimation of the drift $\lambda$: the 
results are valid only in the case when $1/2<H<3/4$. \\

The key point of this method of estimation is that the Hurst exponent $H$ and the diffusion coefficient $\sigma $ can be estimated without
estimating $\lambda$.  We will use this method to estimate the parameters $\sigma$ and $\lambda$. Notice that $H$ was already estimated.\\

Let $\bm{a} = (a_0,\ldots, a_K )$ be a discrete filter of order $L\ge 1$ and length $K+1$, $K \in\IN$ and we require  $L\le K$, i.e.
\[
 \sum_{k=0}^K a_k k^j =0 \quad \hbox{\rm  for } 0\le j\le L-1 \quad\hbox{\rm  and}\quad  \sum_{k=0}^K a_k k^L \ne 0. 
\]
Let it be normalized
\[
 \sum_{k=0}^K (-1)^{1-k} a_k =1.
\]
We will also consider a {\it dilated} filter $\bm{a}^2$ associated to $\bm{a}$. For $0\le k\le K$ we define
\[
a_k^2 = \begin{cases} a_{k'}, & \mbox{if } k=2k' \\0, & \mbox{otherwise.}  \end{cases}
\]
Since $\sum_{k=0}^{2K} a_k^2 k^j=2^j \sum_{k=0}^K a_k k^j $ then the filter $\bm{a}^2$ has the same order than  $\bm{a}$.\\

We are using two filters:
\begin{itemize}
 \item Classical filter. Let $K>0$ and define
\begin{equation*}
  a_k:= \frac{(-1)^{1-k}}{2^k} {K\choose k} =\frac{(-1)^{1-k}}{2^k} \frac{K!}{k!(K-k)!}\qquad \mbox{for }  0\le k\le K.
\end{equation*}

\item Daubechies filters (see \cite{de} for the original definition). The 
 filter is given by 
\begin{equation*}
 \frac{1}{\sqrt{2}} \big(0.48296291314453, -0.8365163037378, 0.22414386804201, 0.12940952255126\big).
\end{equation*}

\end{itemize}

Let $Y^T=(Y_t: 0\le t\le T)$ be  the sample path of the solution of \eqref{mod4}. A discretization of $Y^T$ is
\[
 (X_n:=Y_{n\Delta_N}, n=0,\ldots,N),\qquad 
 N\in \IN,
\]
where $\Delta_N=T/N$ and $N$ is the number of observations of $Y_t$. We denote by
\begin{equation*}
 V_{N,\bm{a}}:= \sum_{i=0}^{N-K}\left( \sum_{k=0}^K a_k X_{i+k} \right)^2,
\end{equation*}
the {\it generalized quadratic variation} associated to the filter $\bm{a}$ (see for instance \cite{is-la}). Then, define the 
following estimators for $H$ and $\sigma$.
\begin{align}
\hat{H}_N&:=\tfrac{1}{2}\log_2\left(\frac{V_{N,\bm{a}^2}}{V_{N,\bm{a}} }\right),\label{est1}\\
\hat{\sigma}_N&:=\left(-2\frac{V_{N,\bm{a}} }{\sum_{k,l} a_ka_l |k-l|^{2 \hat{H}_N } \Delta_N^{2 \hat{H}_N } }\right)^{1/2}.\label{est2}
\end{align}

Brouste and Iacus (see Th. 1 in \cite{br-ia}) have shown the next result. 
\begin{theorem}
Let  be a filter of order $L \ge 2$. Then, both estimators $\hat{H}_N$ and $\hat{\sigma}_N$  are
strongly consistent, i.e.
\[
( \hat{H}_N,\hat{\sigma}_N) \stackrel{\hbox{\rm a.s.}}{\longrightarrow} (H,\sigma) \quad \hbox{\rm as}\quad N \rightarrow+\infty.
\]
Moreover, we have asymptotic normality property:  for all $H \in (0, 1)$,

\begin{align*}
\sqrt{N} (\hat{H}_N  - H ) &\stackrel{\mathcal{L}}{\longrightarrow} N (0, \Gamma_1 (a,\sigma,H)) \quad  \hbox{\rm as} \quad N
\rightarrow +\infty,\\
\frac{\sqrt{N}}{\log N} ( \hat{\sigma}_N  - \sigma ) &\stackrel{\mathcal{L}}{\longrightarrow} N (0, \Gamma_2 (a,\sigma,H)) \quad 
\hbox{\rm as} \quad N \rightarrow +\infty,
\end{align*}

where $\Gamma_1$ and $\Gamma_2$ are symmetric definite positive matrices depending on $\sigma, H$ and the filter $\bm{a}$.
\hfill$\Box$
\end{theorem}

With this result is possible to obtain an estimator for the parameters $\sigma$.

\subsubsection{Estimation of the drift parameter $\lambda$ when both $H$ and $\sigma$ are unknown}

Hu and Nualart, \cite{hu-nu} have shown that 
\[
 \lim_{t\rightarrow \infty} Var (Y_t)= \lim_{t\rightarrow \infty} \frac{1}{t}\int_0^t Y_t^2 dt = \frac{\sigma^2 \Gamma(2H+1) }
 {2\lambda^{2H}}:=\mu_2.
\]
This equation gives a $\lambda$ estimator, namely
\begin{equation}
 \hat{\lambda}_N = \left(\frac{2\hat{\mu}_{2,N}}{\hat{\sigma}_N^2 \Gamma(2\hat{H}_N+1) }  \right)^{-\tfrac{1}{2\hat{H}_N}},\label{est3}
\end{equation}
where $\hat{\mu}_{2,N}$ is the empirical moment of order $2$, i.e
\[
 \hat{\mu}_{2,N} =\tfrac{1}{N}\sum_{n=1}^N X_N^2.
\]
Set $T_N=N\Delta_N$. We have the next result.
\begin{theorem}
 Let $H \in \big(\tfrac{1}{2} , \tfrac{3}{4}\big)$ and a mesh satisfying the condition $N \Delta_N^p\rightarrow 0$, $p>1$,
and $ \Delta_N (log N )^2 \rightarrow 0$ as $N \rightarrow +\infty$. Then, as $N \rightarrow +\infty$,
\[
 \hat{\lambda}_N \stackrel{\hbox{\rm a.s.}}{\longrightarrow}  \lambda,
\]
and 
\[
 \sqrt{T_N} ( \hat{\lambda}_N -\lambda) \stackrel{\mathcal{L}}{\longrightarrow} N (0, \Gamma_3 (\sigma,H)),
\]
where $\Gamma_3 (\sigma,H)=\lambda \big(\tfrac{\sigma_H}{2H} \big)^2$ and 
\[
 \sigma_H^2= (4H+1)\left(1+\frac{\Gamma(1-4H)\Gamma(4H-1) }{\Gamma(2-2H)\Gamma(2H)} \right).
\]
\hfill$\Box$
\end{theorem}
For the proof see Theorem 2 in \cite{br-ia}.\\

\section{Results with the fractional Ornstein-Uhlenbeck model}
\label{re-fou}

In this section we present the estimated mortality rates with the use of the model 
described in section \ref{intro}. We have got the data from the website of Human Mortality Database, we have get the mortality for 
the Italian population between 1950 to 2004.\\

In first place, we present the estimation of the $H$ parameter. In second place, we present the results on simulated mortality rates using  equations \eqref{est2}-\eqref{est3} to estimate the parameters $\sigma,\lambda$. The parameter $\alpha_0$ has been fixed with the use of equation \eqref{alpha0}.\\

We run 10000 simulations of the mortalities rates from ages $0$ to age $90$. To do that, we have simulate a fractional Brownian motion $B_t^{\hat H}$ and using equation \eqref{mod-1} we have estimated the mortality rate. To run the fBm simulations we have used the function {\bf fbm} which is includes in the R library {\bf somebm}.  We also include a  95.5\% confidence interval.\\

We present the results for women and men in sections \ref{re-wom} and \ref{re-men}, respectively. 
 
Estimations and predictions were performed using R Ver. 3.2.3 (R Core Team, 2015), and specialized packages 
Fractal (Time Series Modeling and Analysis
Version 2.0-1, 2016), Pracma (Practical Numerical Math Functions 2.0.7, 2017) and somebm (some Brownian motions simulation function Version 0.1, 2016).

\subsection{Hurst estimation}\label{hu-est}

For estimate the Hurst parameter we have used two R routines: FDWhittle, RoverS and hurstexp. 
The two first routines are from the fractal library, while the latter is from the pracma library.\\

The former routine estimate the Hurst parameter by Whittle's method as was described in the subsection \ref{Desc-Est-H}. RoverS routine
estimate $\hat H $ by rescaled range (R/S) method. The hurstexp routine estimate $\hat H $  using 
R/S analysis. \\

\begin{figure}[H]        
\includegraphics[width = 4.5in,  height=4in]{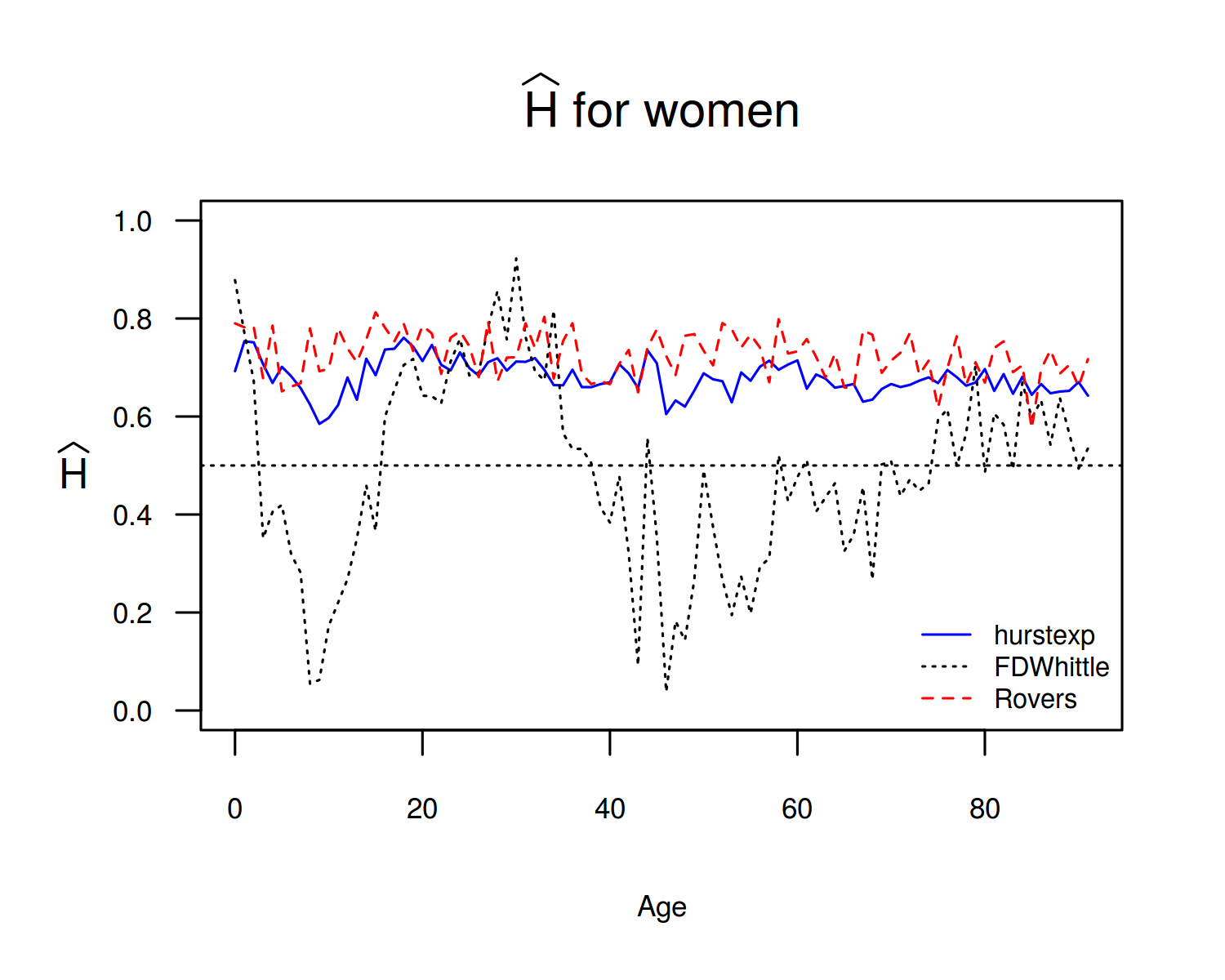}
\caption{\bf Estimated Hurst parameter using R-routines.}
\label{graph-Hurst_Est_Wo}
\end{figure}

\begin{figure}[H]        
\includegraphics[width = 4.5in, height=4in]{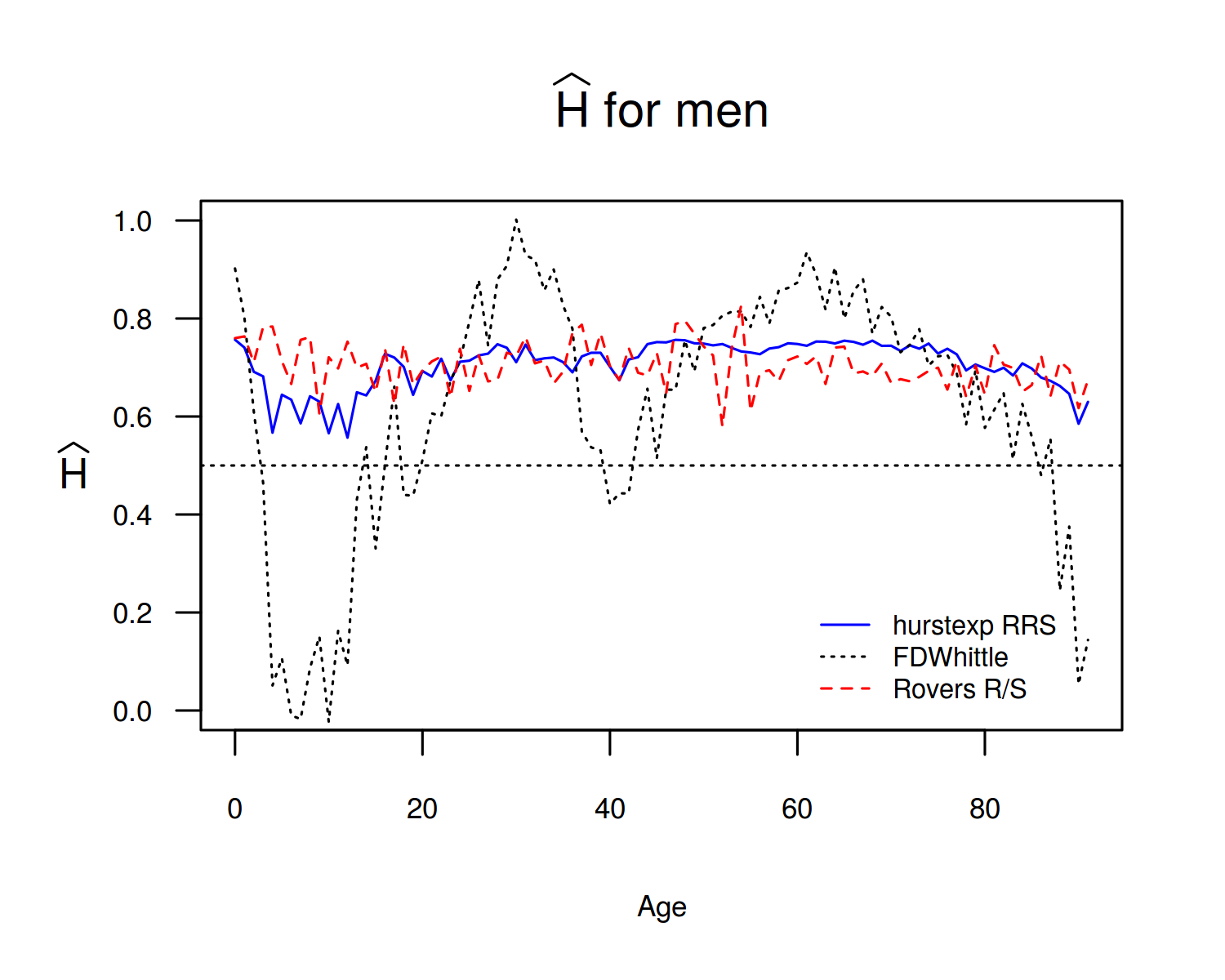}
\caption{\bf Estimated Hurst parameter using R-routines.}
\label{graph-Hurst_Est_Me}
\end{figure}

Finally, figures \ref{graph-Hurst_Est_Wo} and \ref{graph-Hurst_Est_Me} show
 the estimated Hurst parameter for women and men separately.

\vspace*{0.1cm}

With the rescaled range R/S and R/S  methods we obtain a consistent estimator for the Hurst
parameter in the sense that they do not present dramatic changes trough the time. Moreover,
the $H$ estimated with these two methods take values in the interval 
(0.57,0.80) approximately. This tell us that the data has the long memory property as 
was mentioned in section \ref{fgn}. Same results are obtained for the men and women.\\

Notice that the estimated parameters using Whittle method have high variation trough the time, in opposition to those
obtained with the other two methods. So that the estimated Hurst parameters look not good to perform the simulations 
with this method. The high variation on the Hurst estimated values could be explained because Whittle method uses
the periodogram to estimate $H$ while the other methods uses the raw data.\\

Since rescaled range R/S and R/S  methods have estimated very similar $H$, we decide to use the Hurst coefficients obtained with the
method of R/S to perform the mortalities rates simulations. \\

\vspace*{0.1cm}

\subsection{Results for women}\label{re-wom}

\begin{figure}[H]        
\includegraphics[width = 2.85in]{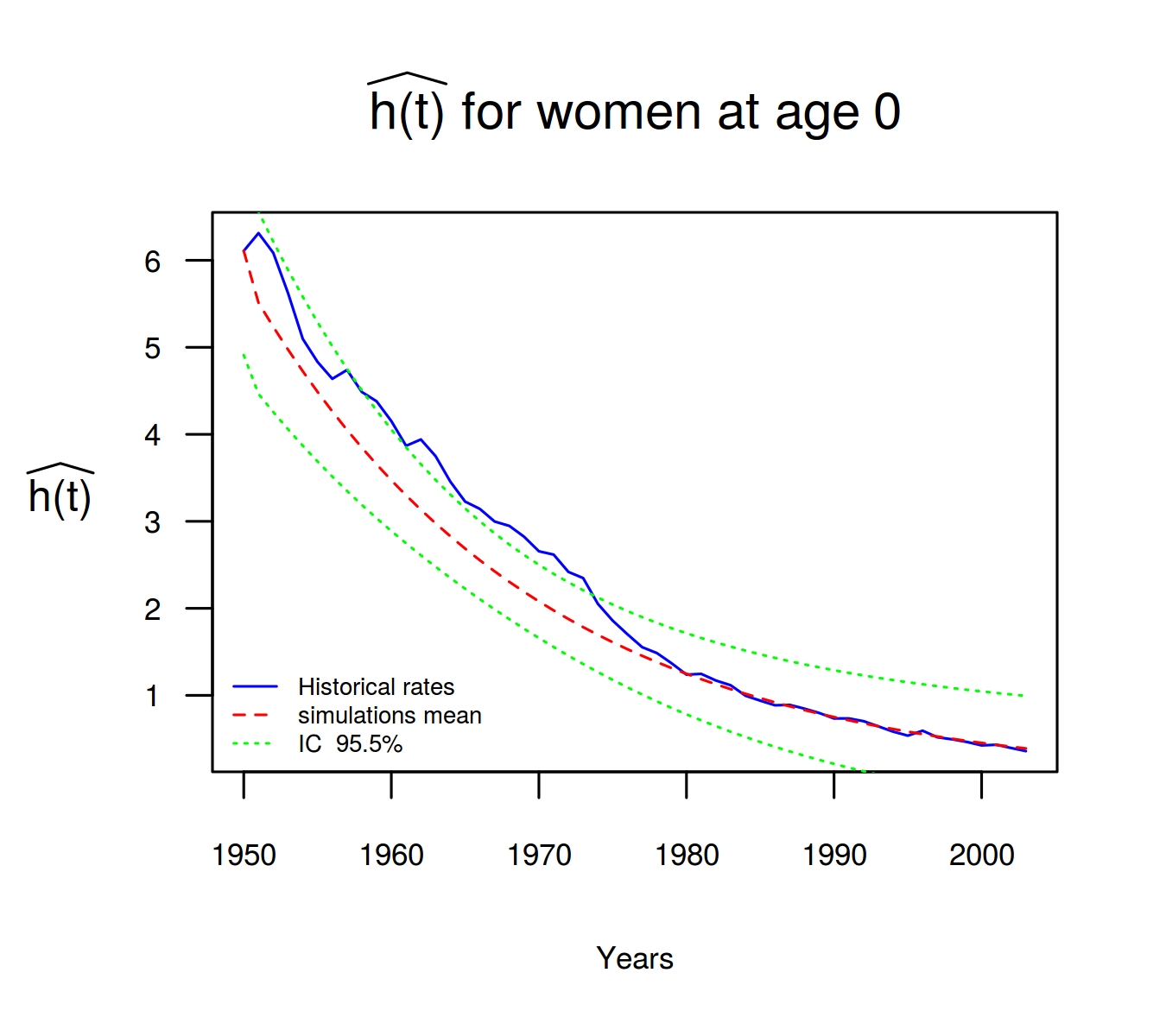}
\includegraphics[width = 2.85in]{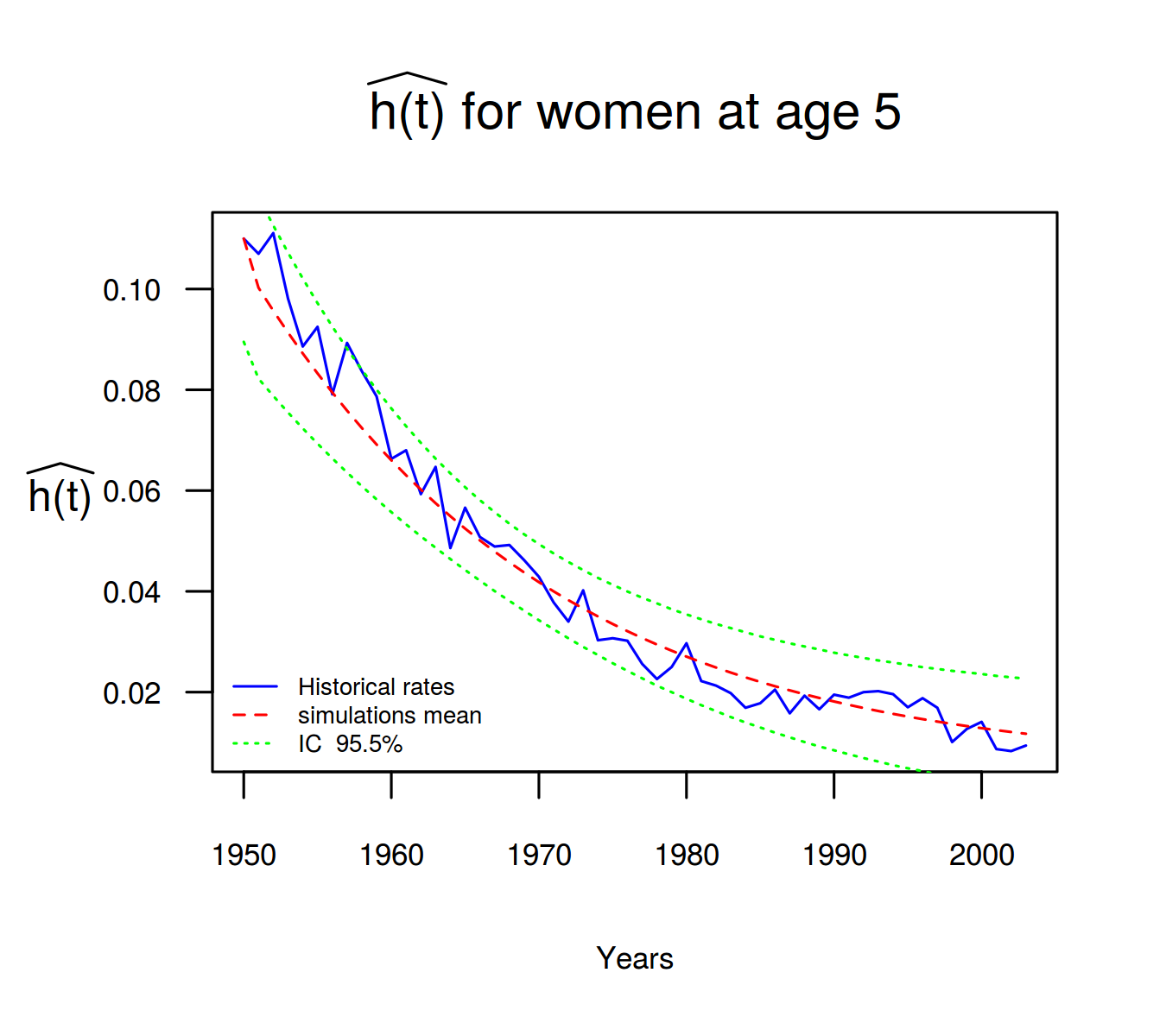} 
\includegraphics[width = 2.85in]{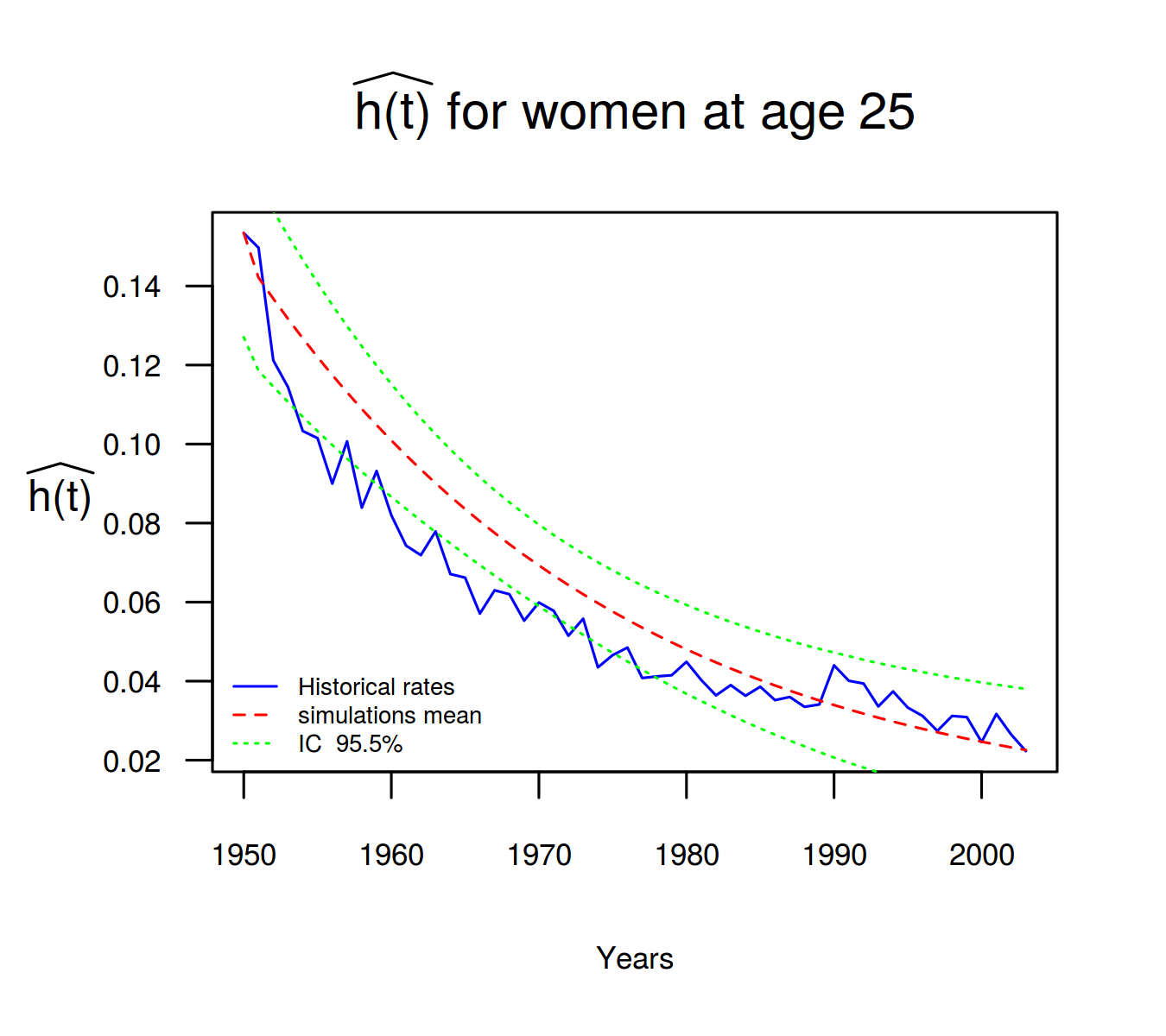} 
\includegraphics[width = 2.85in]{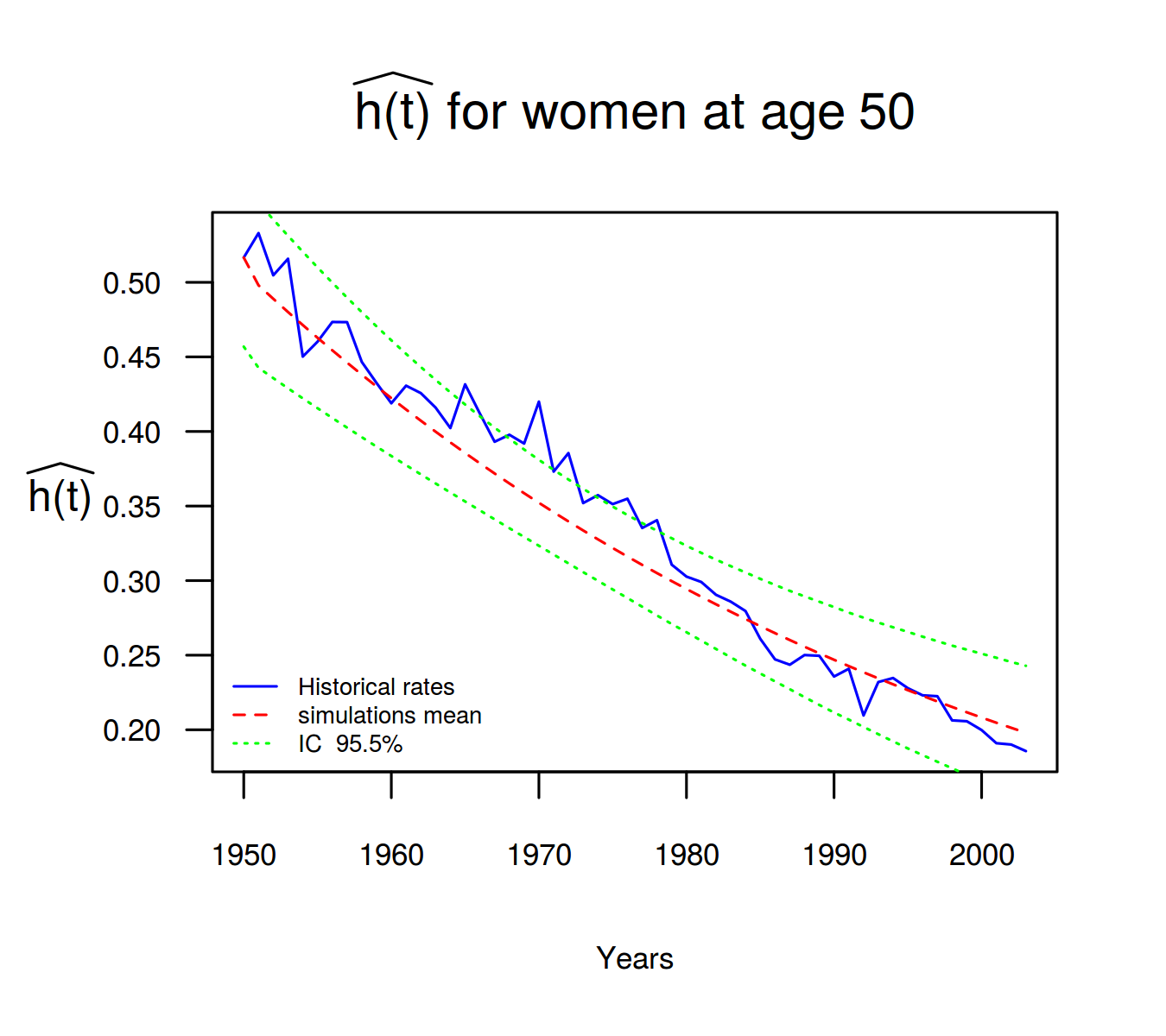} 
\caption{Simulations for the rate mortality with the fOU model: ages $0,5,25,50$ and $N=10000$.}
\label{graph-simu_FOU1}
\end{figure}

\begin{figure}[H]        
\includegraphics[width = 2.85in]{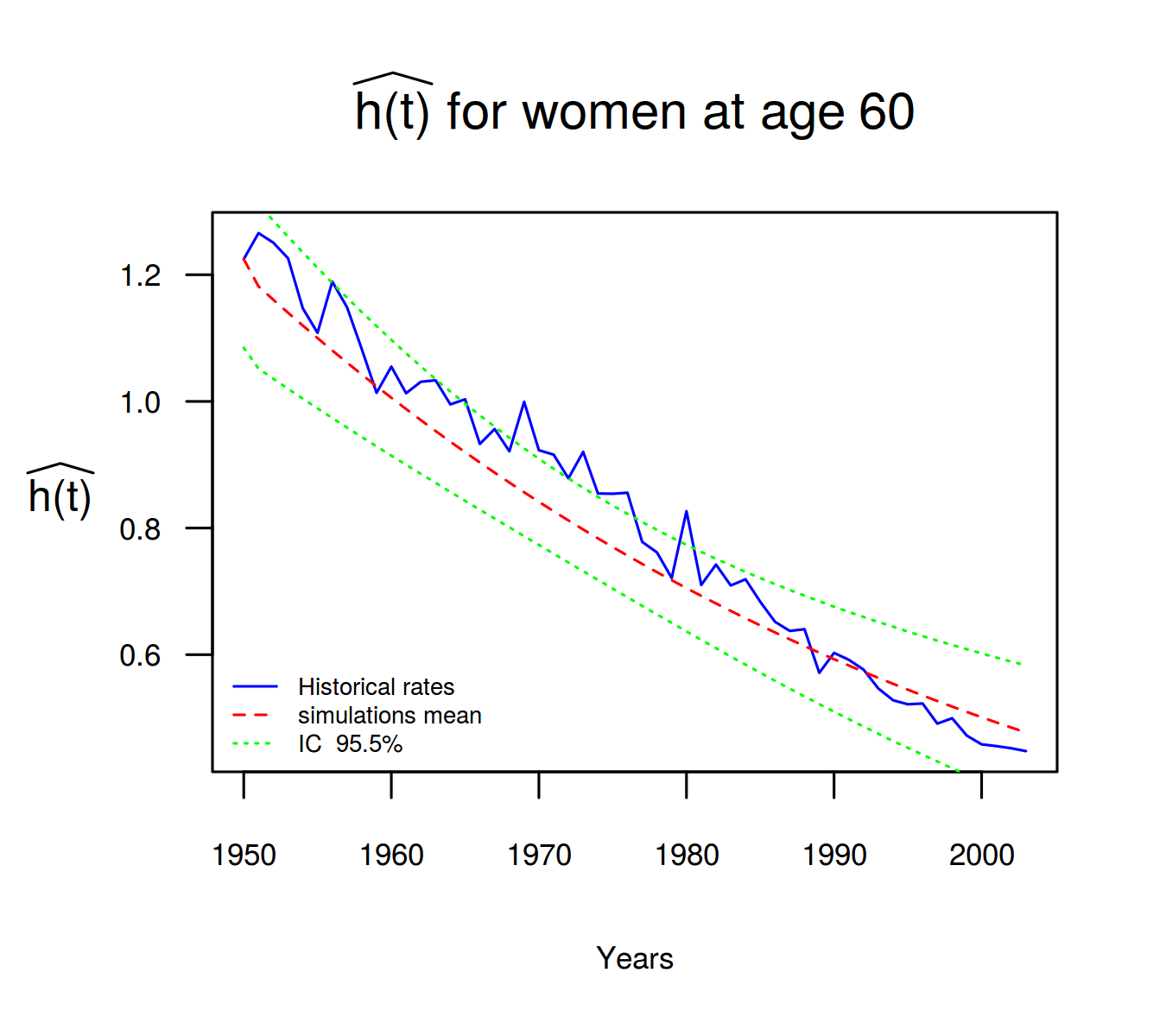}
\includegraphics[width = 2.85in]{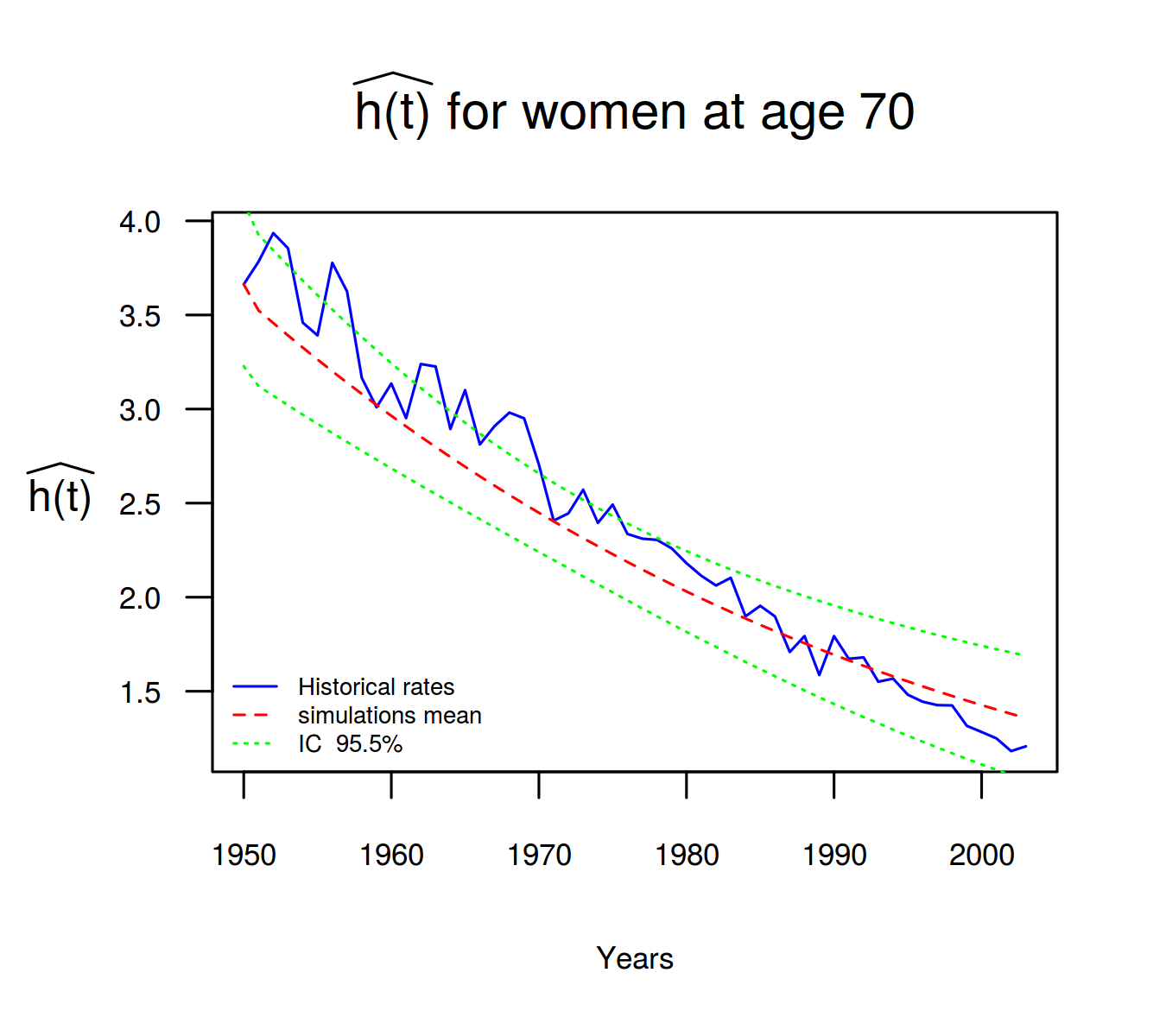}
\includegraphics[width = 2.85in]{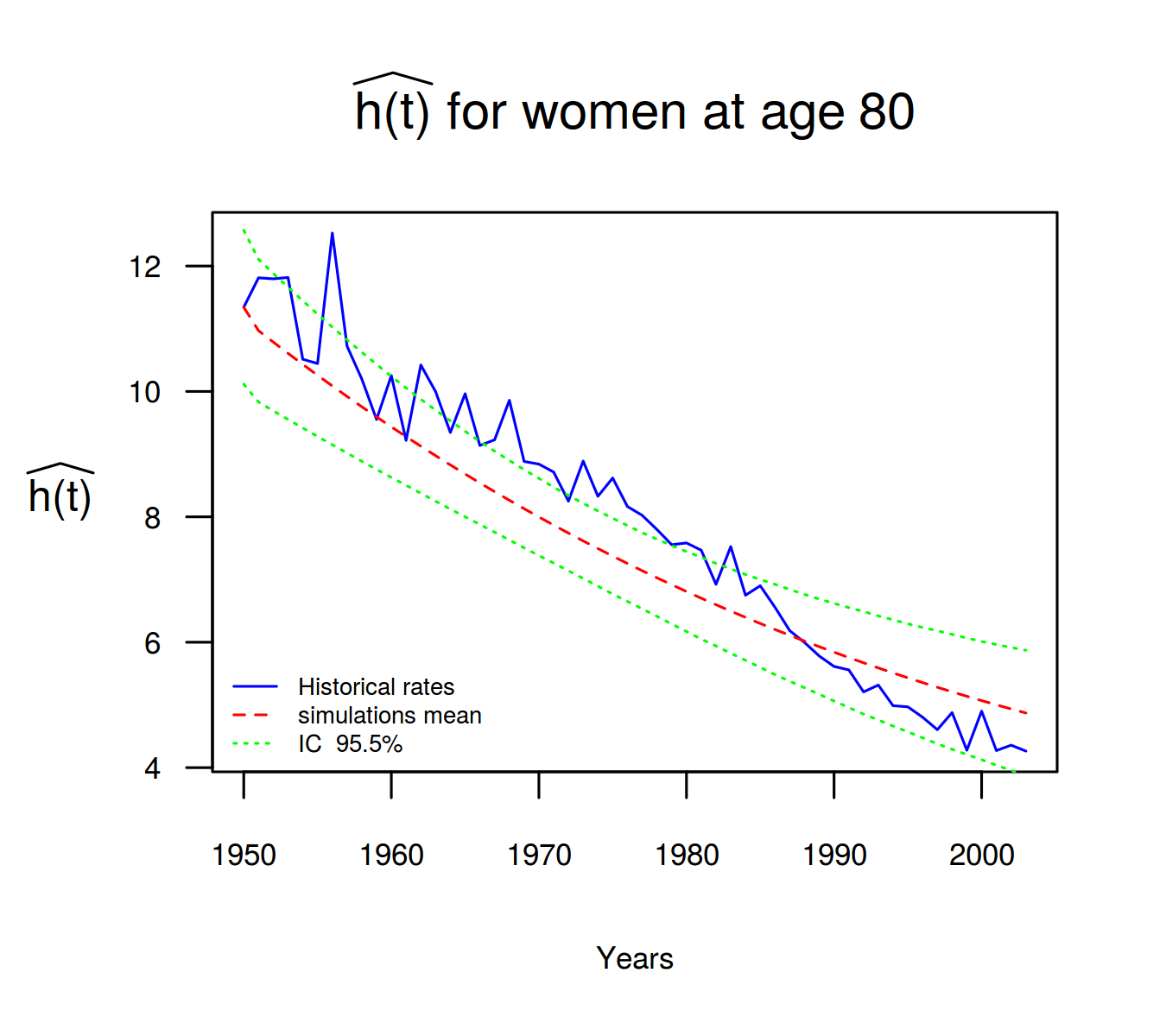}
\includegraphics[width = 2.85in]{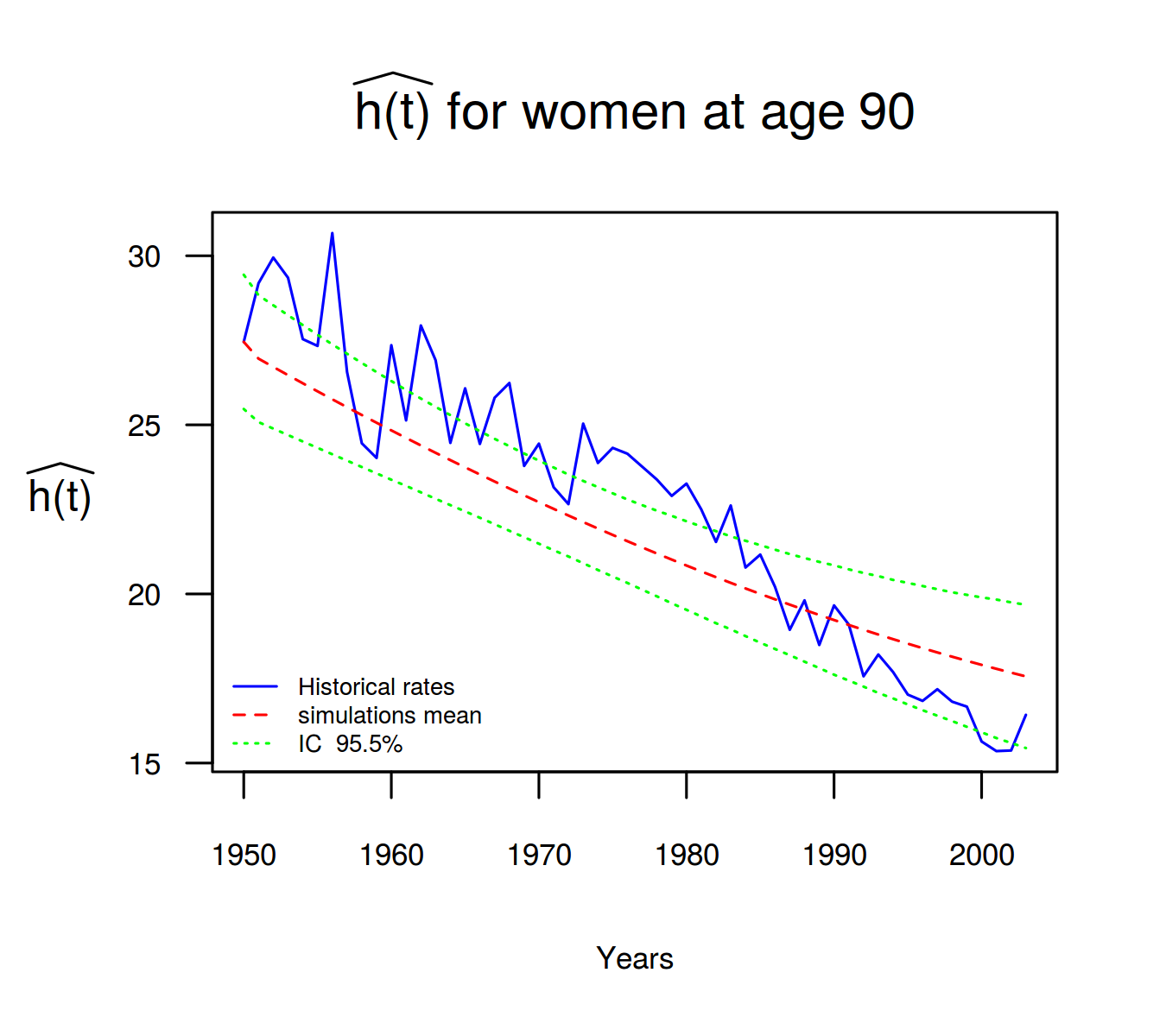}
\caption{\bf Simulations for the rate mortality with the fOU model: ages $60,70,80,90$ and $N=10000$.}
\label{graph-simu_FOU2}
\end{figure}

We present the results for 10000 simulations of the mortalities rates for ages: $0,5,25,50,60,70,80,90$. 
We graph the historical rate mortality, the mean of all simulations and the 95.5\% confidence interval. See figures \ref{graph-simu_FOU1} and  \ref{graph-simu_FOU2}.\\

In general, for all ages, the model is well fitted, in particular, after the 80's. Nevertheless, there are some time periods where the model is not 
so good as we want to.  For instance, for the age 0 and for some ages and between 60's and the 70's years the model underestimates the rate mortality and for age 25 and between late 50's and early 70's overestimates the rate mortality. \\

\vspace*{0.1cm}

From the data we have noticed that, for ages between 20 to 35 approximately (in the beginning of the 90's) the phenomena of AIDS has made a little increase in the rates mortality; however, the model captures this circumstance well, since the estimation still remains inside 
the confidence intervals. We believed that in order  to compensate this variation, the model overestimates mortalities rates 
 for the period between late 50's and first 70's. \\

For older ages (see figure \ref{graph-simu_FOU2}) we observe that for ages 60 and 70 the estimation is well fitted trough the years. we notice that 
the predicted rates is not so far away and that the historical rates are inside the confidence interval. For the very oldest ages the estimation is not so good as for earlier ages. The main difference is in 50's when the absolute number of living persons arriving to those ages were small. \\

All this suggest that a better model could include a short and a long-term memory process, so that the model could help us to control the short-term variations in a better way.

\vspace*{0.1cm}

 \subsection{Results for men}\label{re-men}

\begin{figure}[H]        
\includegraphics[width = 2.85in]{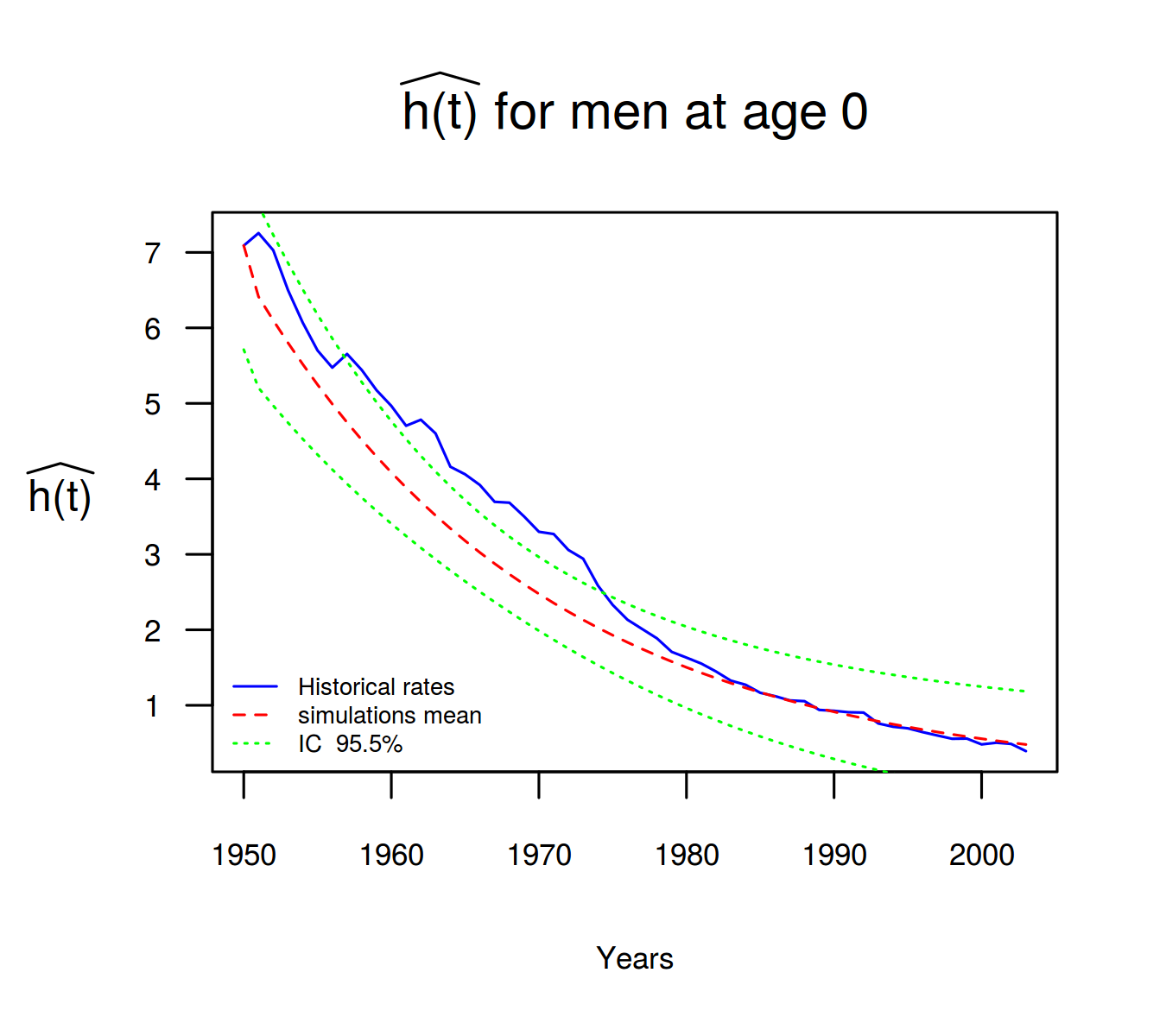}
\includegraphics[width = 2.85in]{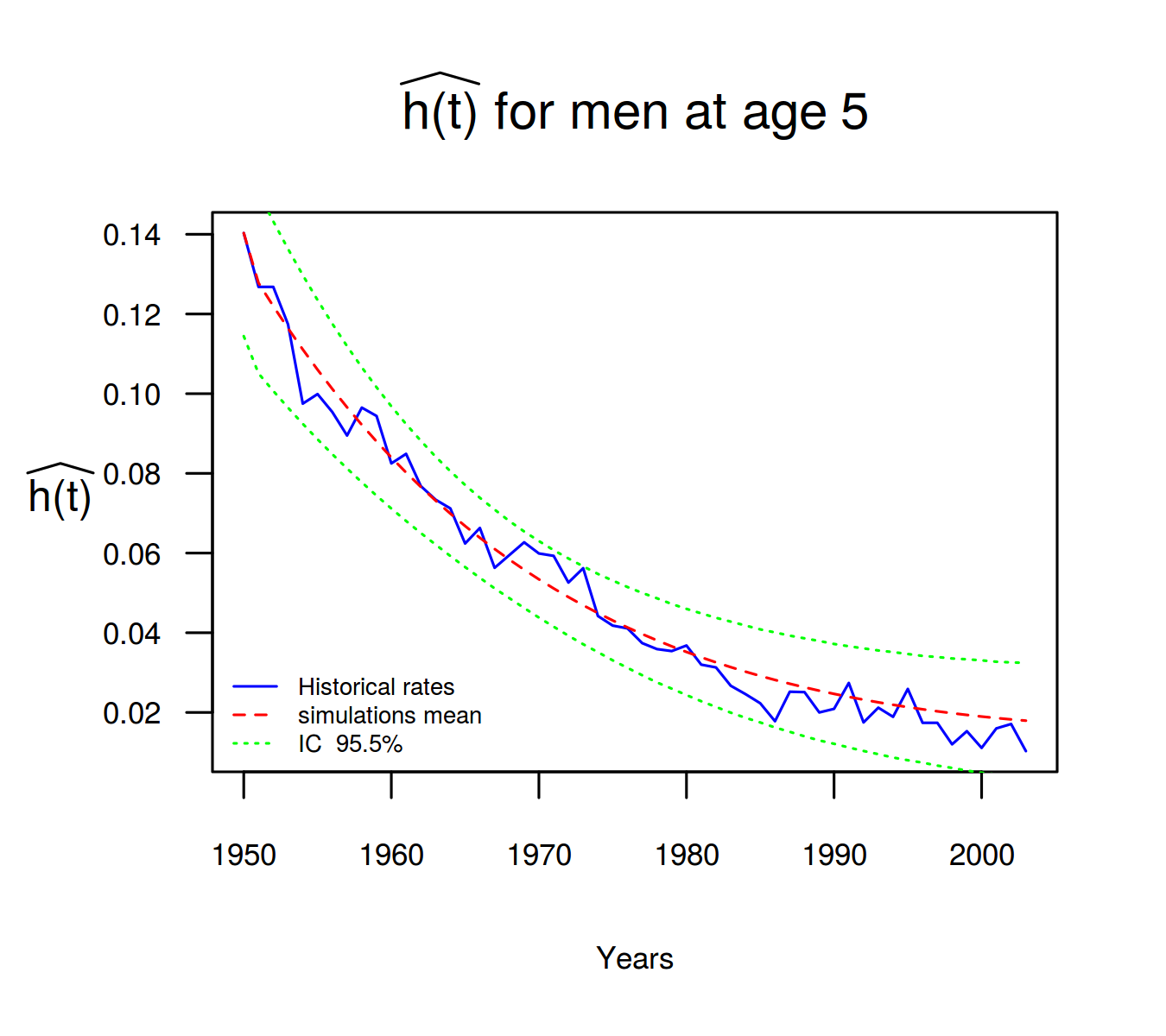} 
\includegraphics[width = 2.85in]{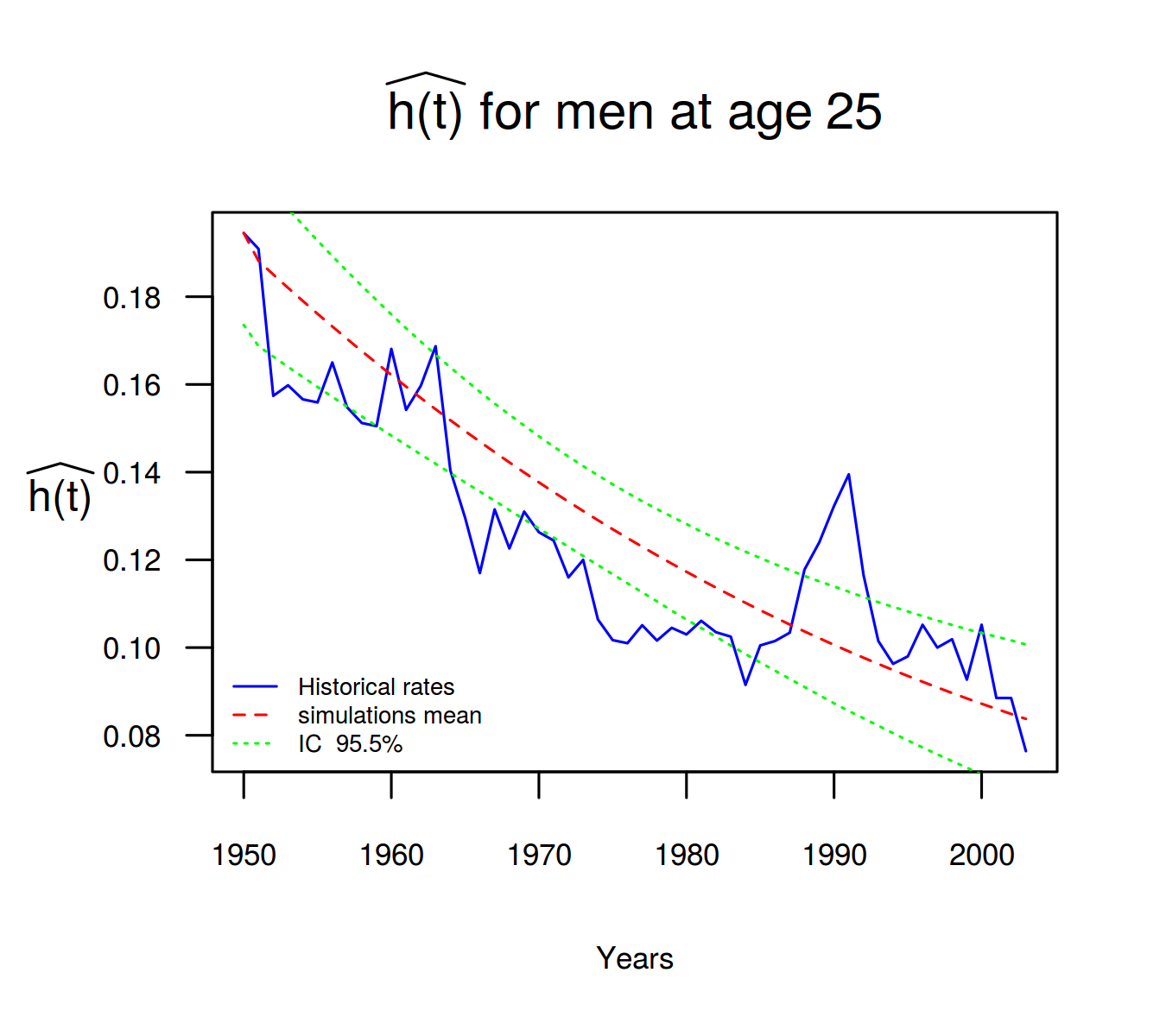} 
\includegraphics[width = 2.85in]{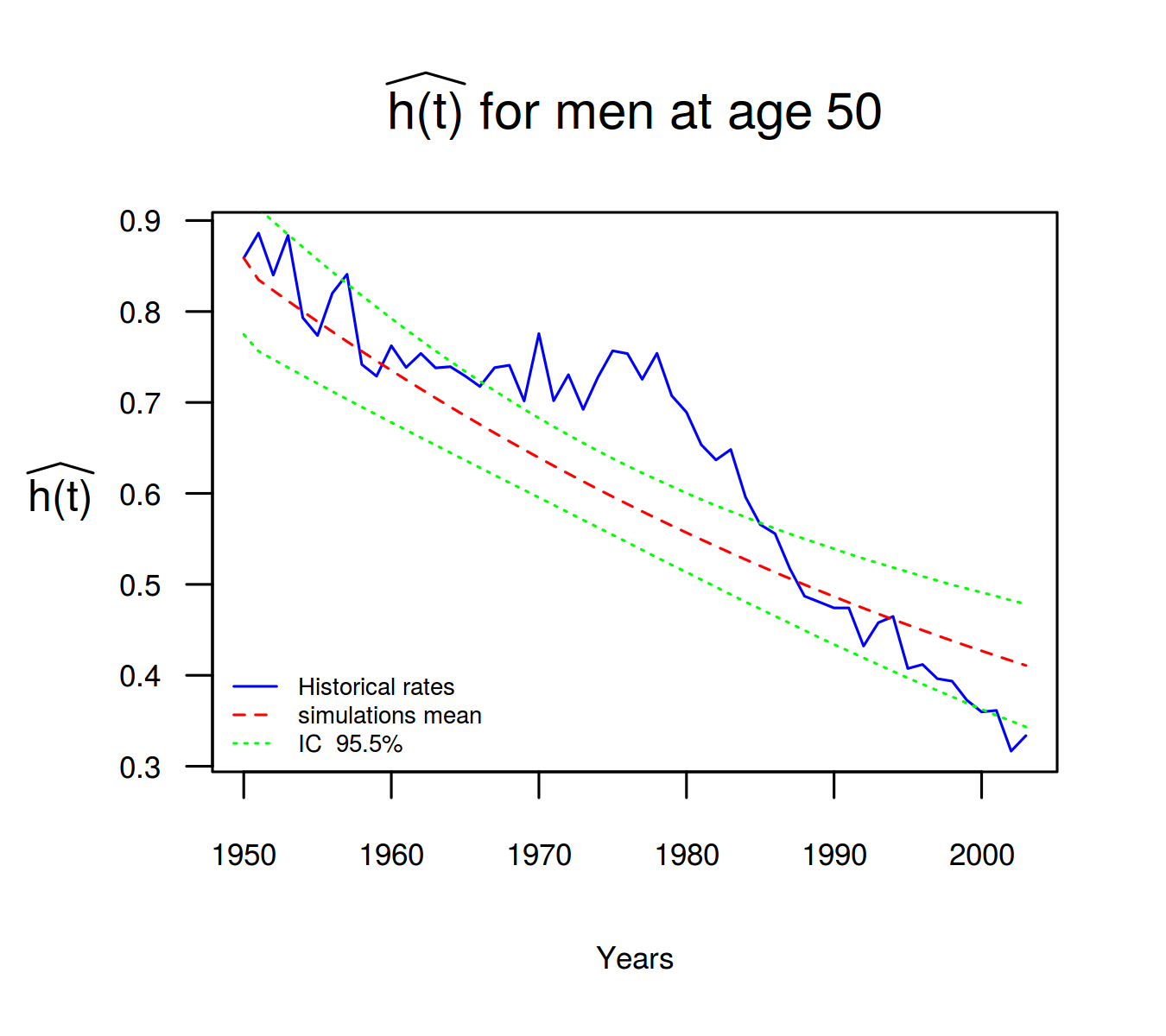} 
\caption{\bf Simulations for the rate mortality with the fOU model: ages $0,5$ and $N=10000$.}
\label{graph-simu_FOU3}
\end{figure}

\begin{figure}[H]        
\includegraphics[width = 2.85in]{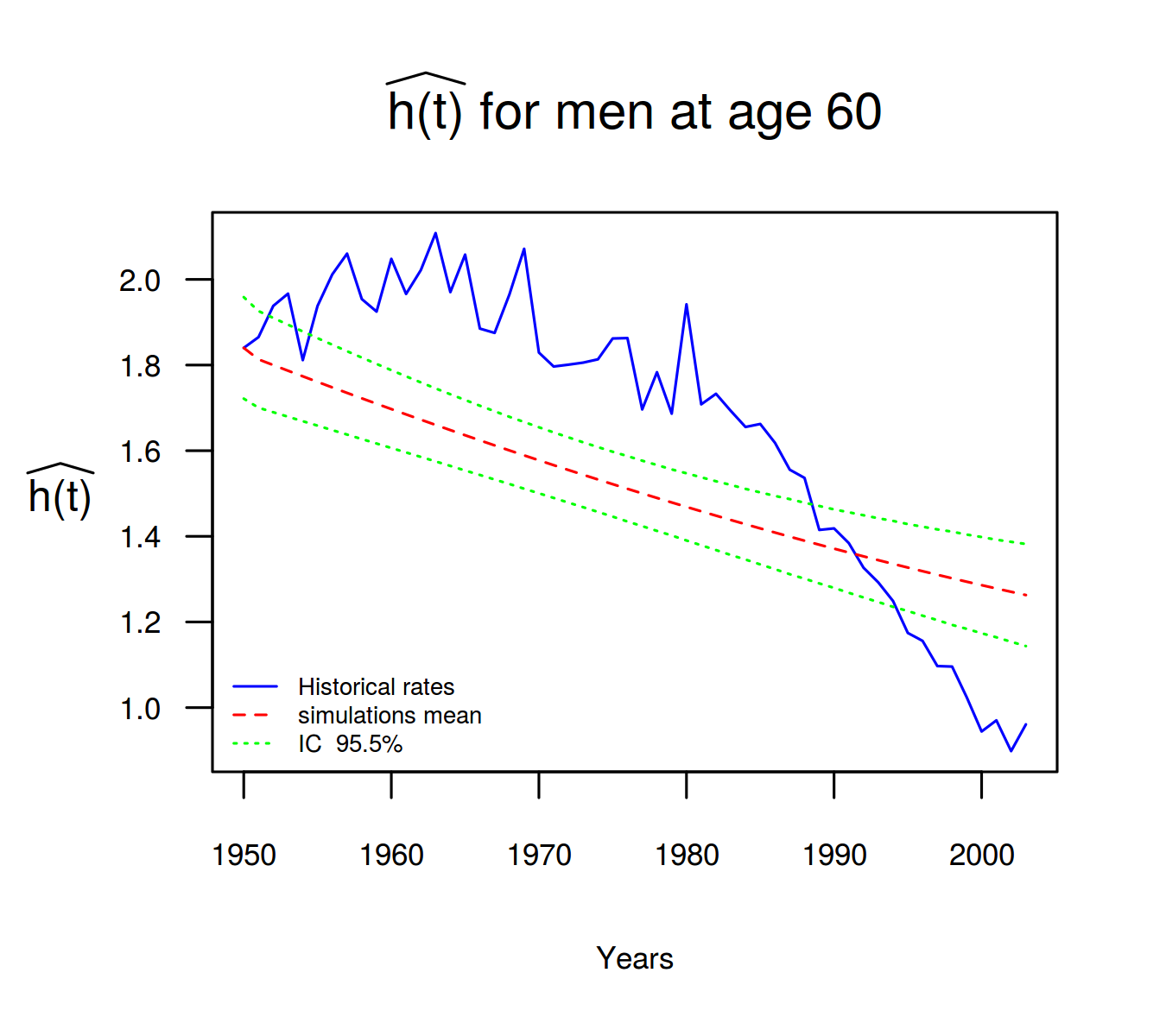}
\includegraphics[width = 2.85in]{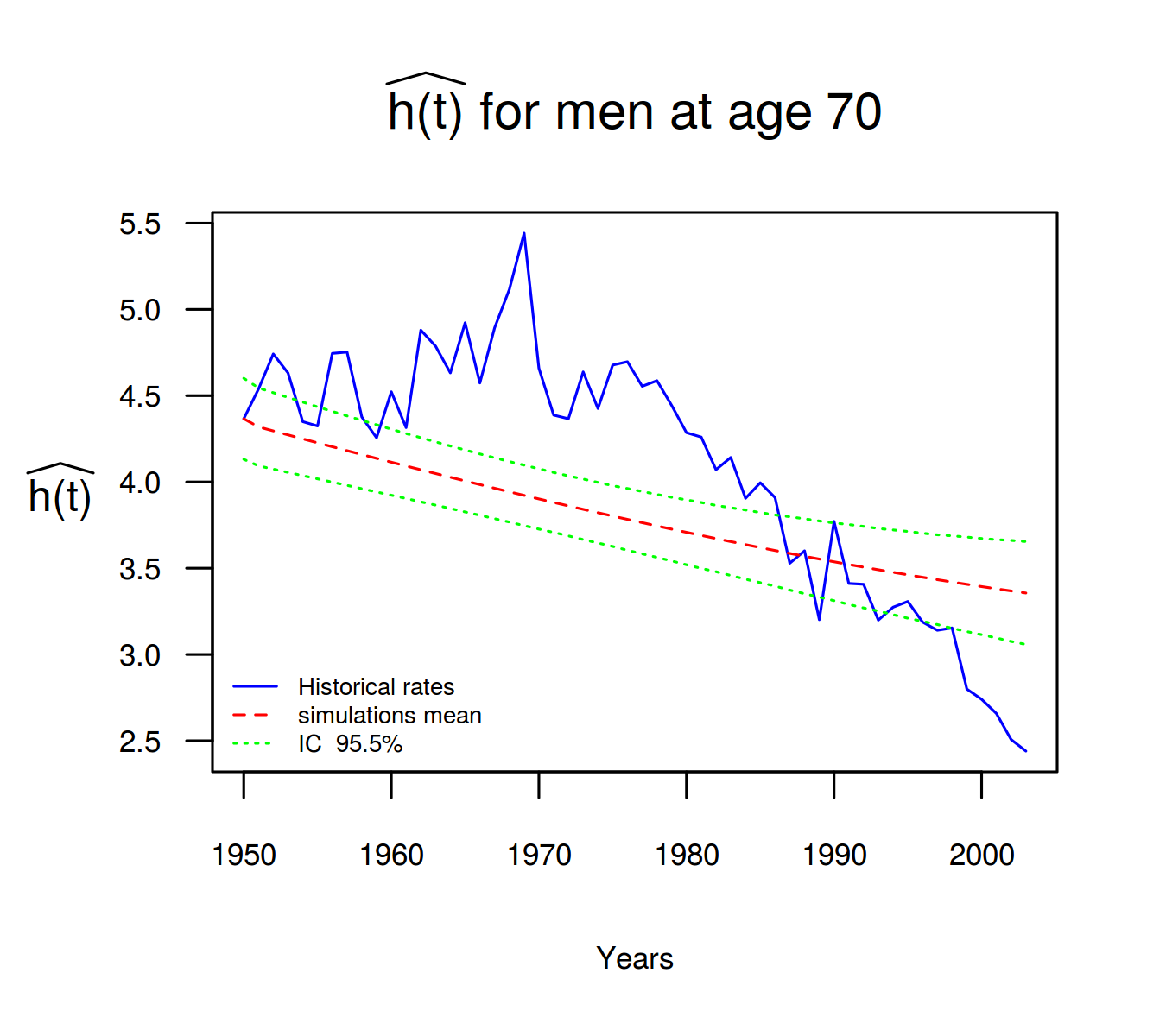}
\includegraphics[width = 2.85in]{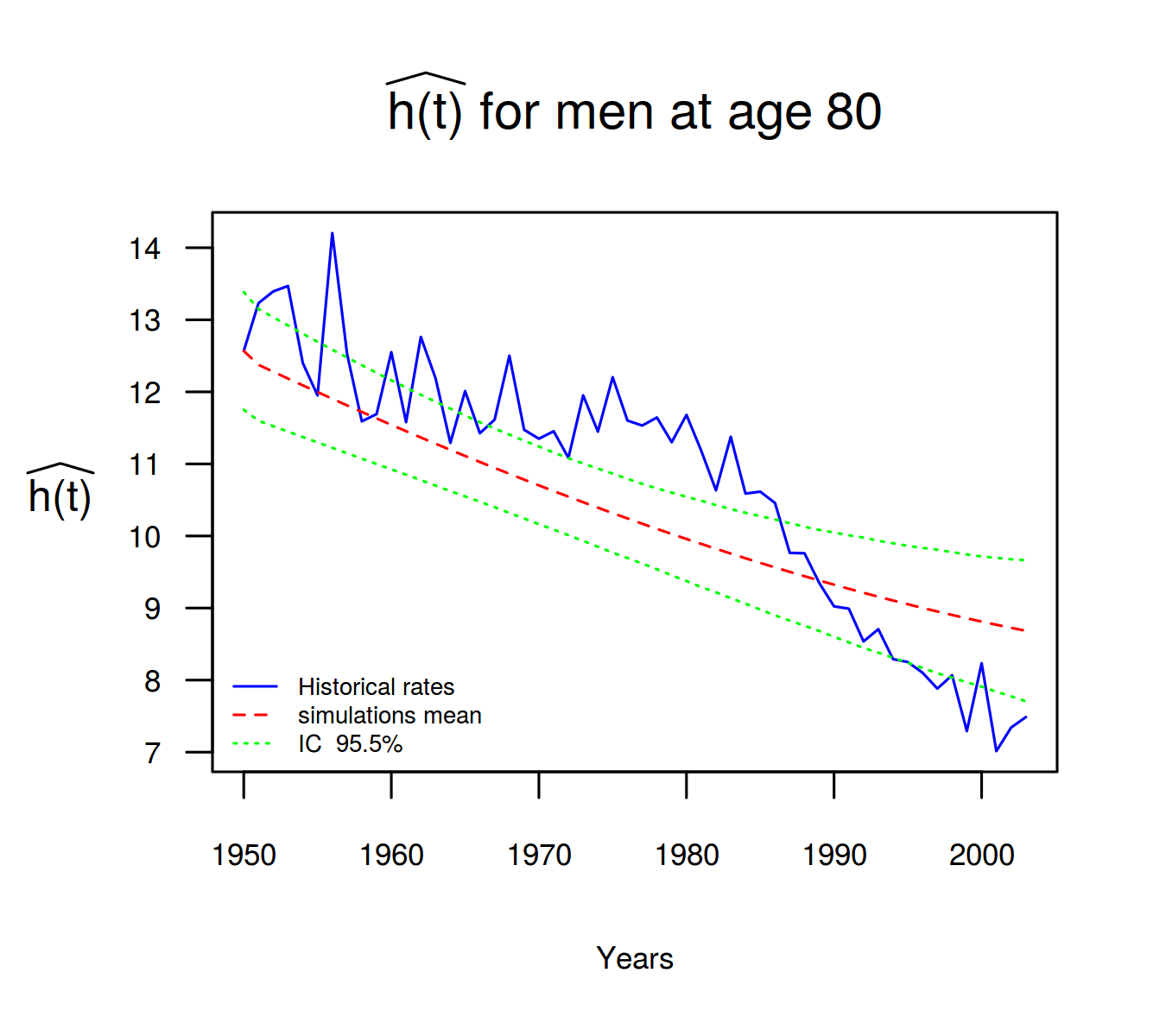}
\includegraphics[width = 2.85in]{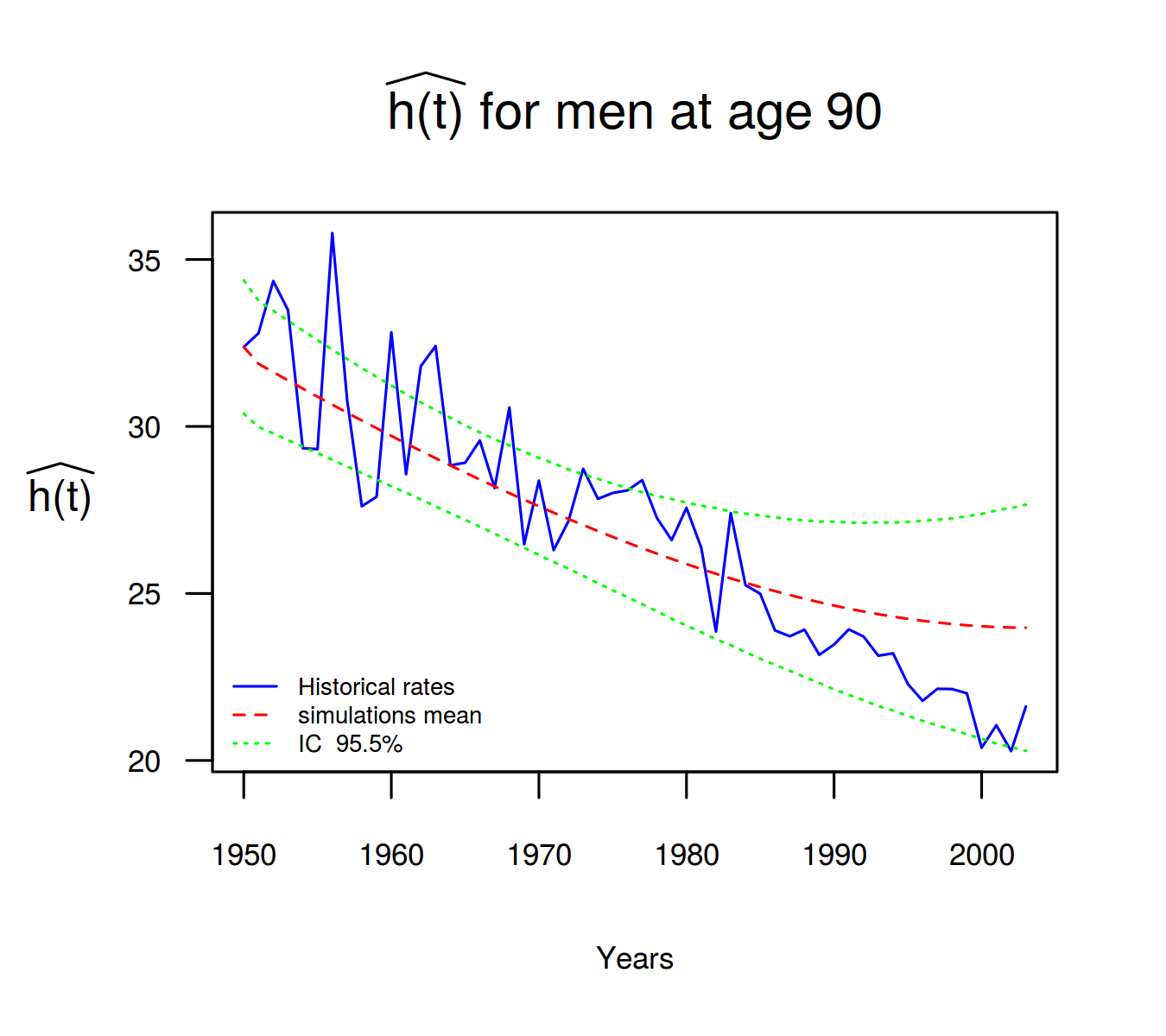}
\caption{\bf Simulations for the rate mortality with the fOU model: ages $60,70,80,90$ and $N=10000$.}
\label{graph-simu_FOU4}
\end{figure}

As in the case for women, we present results for 10000 simulations of the mortalities rates for ages: $0,5,25,50,60,70,80,90$. 
We graph the historical rate mortality, the mean of all simulations and the 95.5\% confidence interval. See figures \ref{graph-simu_FOU3} and  \ref{graph-simu_FOU4}.\\

\vspace*{0.1cm}

As in the case for women, the proposed model for men is well fitted. The AIDS pandemic is also noticed for ages between 25 to 35 and
it generates an increase in the rates mortalities for these ages; As a matter fact, this increase is heavier
than the one for the women, this has caused a overestimation in the first 35 years and latter a underestimation 
of the rates mortalities. As was mentioned before, if we include in the model a short-term process, the model could be better fitted.

\vspace*{0.1cm}

For older ages, we observed that the data is irregular, so it is necessary to use a more complex model to 
fit this data.

\section{Conclusions}

We have adjust a model for the Italian mortality rates with a geometric-type fractional Ornstein-Uhlenbeck process. Our 
main hypothesis was that, for a fixed age, the mortality rates changes trough the time slowly, so that a stochastic differential equations
that captures the long-range dependence could be a good model. We have use a stochastic differential equation
with a fractional Brownian motion as a driven noise with $H\in (0.5,1)$ in order to satisfy the long-range dependence property. With the data we have fixed the Hurst coefficient and we have confirmed our hypothesis since 
we have found that the estimated Hurst is in $(0.58,0.8)$.  \\

Notice that we have consider a more general model that the one used in \cite{gi-or-be}, this is because we have included the possibility that the 
Hurst parameter could be equal to $1/2$, which is the case when the fractional Brownian motion becomes a standard Brownian motion. Therefore, 
when $H=1/2$ we recover the  Giacometti, Ortobelli and Bertocchi model.\\


The model is, in general, well behaved. For some ages we found some shortcomings that suggest that the use of more terms could improve
the model. The long-range dependence model proposed in this paper is good enough to reproduce the mortality rates. If we add some extra terms to make it more flexible to reproduce the cases where the mortality rates have more variations then it will generates a more accurate model. We are starting to work on this extension of the model. Moreover, a multiplicative noise model will be subject of a future research.\\

\end{document}